\documentclass{amsart}

\usepackage{amsthm,amsfonts,amsmath,amssymb,mathrsfs,verbatim,cite}

\numberwithin{equation}{section}

\newtheorem{lemma}{Lemma}[section]
\newtheorem{theorem}{Theorem}[section]
\newtheorem{corollary}{Corollary}[section]

\newtheorem{definition}{Definition}[section]

\renewcommand{\Re}{\textup{Re}}
\newcommand{\dup}{\textup{d}}
\newcommand{\eup}{\textup{e}}
\newcommand{\iup}{\textup{i}}
\newcommand{\Int}{\int\limits}
\newcommand{\abs}[1]{\lvert#1\rvert}
\newcommand{\Abs}[1]{\bigl\lvert#1\bigr\rvert}
\newcommand{\F}{\mathbb F}
\newcommand{\Q}{\mathbb Q}
\newcommand{\R}{\mathbb R}
\newcommand{\Z}{\mathbb Z}
\newcommand{\Complex}{\mathbb C}
\newcommand{\A}[1]{\textup{A}_{#1}}
\newcommand{\la}{\lambda}
\newcommand{\La}{\Lambda}
\newcommand{\las}{\lambda^{(s)}}
\newcommand{\laa}[1]{\lambda^{(#1)}}
\newcommand{\mus}{\mu^{(s)}}
\newcommand{\muu}[1]{\mu^{(#1)}}
\newcommand{\xs}{x^{(s)}}
\newcommand{\ba}{\bar\alpha}
\newcommand{\xa}[1]{x^{(#1)}}
\newcommand{\gsl}{\mathfrak{sl}}
\newcommand{\g}{\mathfrak{g}}
\newcommand{\h}{\mathfrak{h}}
\newcommand{\Symm}{\mathfrak{S}}
\newcommand{\bil}[2]{(#1,#2)}

\begin{document}

\title{A Selberg integral for the Lie algebra $\A{n}$}

\author{S. Ole Warnaar}\thanks{Work supported by the Australian Research Council}
\address{Department of Mathematics and Statistics,
The University of Melbourne, VIC 3010, Australia}
\email{warnaar@ms.unimelb.edu.au}



\begin{abstract}
A new $q$-binomial theorem for Macdonald polynomials is employed 
to prove an $\A{n}$ analogue of the celebrated Selberg integral.
This confirms the $\g=\A{n}$ case of a conjecture by Mukhin and Varchenko 
concerning the existence of a Selberg integral for every simple Lie
algebra $\g$.
\end{abstract}

\maketitle

\section{Introduction}

\subsection{$\g$-Selberg integrals}
In 1944 Selberg published the following remarkable
multiple integral \cite{Selberg44}.
Let $k$ be a positive integer, $t=(t_1,\dots,t_k)$, 
$\dup t=\dup t_1\cdots\dup t_k$,
and 
\begin{equation*}
\Delta(t)=\prod_{1\leq i<j\leq k}(t_i-t_j)
\end{equation*}
the Vandermonde product.
\begin{theorem}[Selberg integral]\label{thmSel}
For $\alpha,\beta,\gamma\in\Complex$ such that 
\begin{equation*}
\Re(\alpha)>0,~\Re(\beta)>0,~\Re(\gamma)>
-\min\{1/k,\Re(\alpha)/(k-1),\Re(\beta)/(k-1)\}
\end{equation*}
there holds
\begin{multline}\label{Selberg}
\Int_{[0,1]^k} \abs{\Delta(t)}^{2\gamma}
\prod_{i=1}^k t_i^{\alpha-1} (1-t_i)^{\beta-1}\, \dup t  \\
=\prod_{i=1}^k
\frac{\Gamma(\alpha+(i-1)\gamma) \Gamma(\beta+(i-1)\gamma)\Gamma(i\gamma+1)}
{\Gamma(\alpha+\beta+(i+k-2)\gamma)\Gamma(\gamma+1)}.
\end{multline}
\end{theorem}
When $k=1$ the Selberg integral simplifies to the
Euler beta integral \cite{Euler1730}
\begin{equation}\label{Euler}
\int_0^1 t^{\alpha-1} (1-t)^{\beta-1}\, \dup t =
\frac{\Gamma(\alpha) \Gamma(\beta)}{\Gamma(\alpha+\beta)},
\qquad \Re(\alpha)>0,~\Re(\beta)>0,
\end{equation}
which reduces to the standard definition of the gamma function 
\begin{equation*}
\Gamma(\alpha)=\int_0^{\infty}t^{\alpha-1}\eup^{-t} \dup t,
\qquad \Re(\alpha)>0
\end{equation*}
upon taking $(\beta,t)\to(\zeta,t/\zeta)$ (with $\zeta\in\R$) and letting
$\zeta\to\infty$.

At the time of its publication the Selberg integral was largely overlooked,
but now, more than $60$ years  later, it is widely regarded as
one of the most fundamental and important hypergeometric integrals.
It has connections and applications to orthogonal polynomials,
random matrices, finite reflection groups, 
hyperplane arrangements, Knizhnik--Zamolodchikov equations
and more,
see e.g., \cite{AAR99,Askey80,DX01,EFK03,Macdonald82,Mehta04,Varchenko03,Varchenko04}.

Because of the appearance of the Vandermonde product, the Selberg 
integral may be associated with the root system $\A{k-1}$.
That such a viewpoint is useful is evidenced by Macdonald's famous
ex-conjecture, which attaches a Selberg integral
to any finite reflextion group $G$ \cite{Macdonald82}.
To be precise, Macdonald conjectured a generalisation to $G$
of the exponential limit of Theorem~\ref{thmSel}, known as Mehta's integral: 
\[
\frac{1}{(2\pi)^{k/2}}
\Int_{\R^k} \abs{\Delta(t)}^{2\gamma}
\eup^{-\frac{1}{2}\sum_{i=1}^n t_i^2}\, \dup t  \\
=\prod_{i=1}^k \frac{\Gamma(i\gamma+1)}{\Gamma(\gamma+1)},
\]
see also \cite{Garvan89,Macdonald82,Mehta04,Opdam89,Opdam93}.

\medskip

A different point of view --- and one we wish to adopt in this paper ---
arises from the intimate connection between Knizhnik--Zamolodchikov (KZ) 
equations and hypergeometric integrals \cite{EFK03,SV91,MV00}.
Let $\g$ be a simple Lie algebra of rank $n$, with simple roots and
Chevalley generators given by $\ba_i$ and $e_i,f_i,h_i$ 
for $1\leq i\leq n$.\footnote{We use $\ba_i$ instead of 
the usual $\alpha_i$ to denote the 
simple roots to avoid a clash of notation with the exponents $\alpha_i$ in the
$\A{n}$ Selberg integral of Theorem~\ref{thmSelbergAn}.} 
Let $V_{\la}$ and $V_{\mu}$ be highest weight representations of $\g$
with highest weights $\la$ and $\mu$, and let $u=u(z,w)$ be 
a function with values in $V_{\la}\otimes V_{\mu}$ solving the
KZ equation
\[
\kappa\,\frac{\partial u}{\partial z}=\frac{\Omega}{z-w}\, u,
\qquad\quad
\kappa\,\frac{\partial u}{\partial w}=\frac{\Omega}{w-z}\, u,
\]
where $\Omega$ is the Casimir element. Solutions $u$ with values in 
the space of singular vectors of weight $\la+\mu-\sum_{i=1}^n k_i \ba_i$
are expressible in terms of $k:=k_1+\cdots+k_n$ dimensional integrals of 
hypergeometric type as follows \cite{SV91}:
\[
u(z,w)=\sum u_{I,J}(z,w) \, f^I v_{\la}\otimes f^J v_{\mu}
\]
with
\[
u_{I,J}(z,w)=\int_{\gamma}\Psi(z,w;t)\omega_{I,J}(z,w;t)\dup t.
\]
In the above the sum is over all ordered multisets $I$ and $J$ with elements
taken from $\{1,\dots,n\}$ such that their union contains the
number $i$ exactly $k_i$ times, $v_{\la}$ and $v_{\mu}$ are the 
highest weight vectors of $V_{\la}$ and $V_{\mu}$,
$f^{I} v=(\prod_{i\in I} f_i) v$, $t=(t_1,\dots,t_k)$, 
$\dup t=\dup t_1\cdots \dup t_k$ and $\gamma$ is a suitable integration cycle. 
The function $\omega_{I,J}$ is a rational function that will
not concern us here and $\Psi$, known as the phase function, is
defined as follows. The first $k_1$ integration variables are
attached to the simple root $\ba_1$, the next $k_2$ integration variables
are attached to the simple root $\ba_2$, and so on, such that
$\ba_{t_j}:=\ba_i$ if $k_{i-1}<j\leq k_i$. Then
\begin{multline*}
\Psi(z,w;t)=(z-w)^{\bil{\la}{\mu}/\kappa}
\prod_{i=1}^k (t_i-z)^{-\bil{\la}{\ba_{t_i}}/\kappa}
(t_i-w)^{-\bil{\mu}{\ba_{t_i}}/\kappa} \\
\times
\prod_{1\leq i<j\leq n}(t_i-t_j)^{\bil{\ba_{t_i}}{\ba_{t_j}}/\kappa},
\end{multline*}
with $\bil{\,}{}$ the bilinear symmetric form on $\h^{\ast}$
(the space dual to the Cartan subalgebra $\h$)
normalised such that $\bil{\theta}{\theta}=2$ for the
maximal root $\theta$.

In \cite{MV00} Mukhin and Varchenko formulated a remarkable
conjecture regarding the normalised phase function
\[
\Psi(t)=
\prod_{i=1}^k t_i^{-\bil{\la}{\ba_{t_i}}/\kappa}
(1-t_i)^{-\bil{\mu}{\ba_{t_i}}/\kappa}
\prod_{1\leq i<j\leq n}\abs{t_i-t_j}^{\bil{\ba_{t_i}}{\ba_{t_j}}/\kappa}.
\]
They proposed that if the space of singular vectors of weight
$\la+\mu-\sum_{i=1}^n k_i \ba_i$ is one-dimensional, then
the integral
\begin{equation}\label{MV}
\int_{\Delta} \Psi(t) \dup t
\end{equation}
(with $\Delta\subset[0,1]^k$ an appropriate integration domain
not explicitly given) is expressible as a product of gamma functions.
The original Selberg integral corresponds to the case $\g=\gsl_2$ of 
the Mukhin--Varchenko conjecture.

In the following we restrict our attention 
to $\g=\gsl_{n+1}=\A{n}$, with fundamental weights $\La_1,\dots,\La_n$; 
$\bil{\La_i}{\ba_j}=\delta_{ij}$.
If the weights of $V_{\la}$ and $V_{\mu}$ are
$\la=\sum_{i=1}^n \la_i \La_i$ and 
$\mu=\sum_{i=1}^n \mu_i \La_i$, and if we write
$t=(t_1,\dots,t_k)$ as 
$t=(t^{(1)},\dots,t^{(n)})$, with $t^{(s)}=(t_1^{(s)},\dots,t_{k_s}^{(s)})$
the variables attached to the simple root $\ba_s$, then
\[
\Psi(t)=
\prod_{s=1}^n\biggl[ \Abs{\Delta\bigl(t^{(s)}\bigr)}^{2/\kappa}
\prod_{i=1}^{k_s}
\bigl(t_i^{(s)}\bigr)^{-\la_s/\kappa}
\bigl(1-t_i^{(s)}\bigr)^{-\mu_s/\kappa} \biggr]
\prod_{s=1}^{n-1} 
\Abs{\Delta\bigl(t^{(s)},t^{(s+1)}\bigr)}^{-1/\kappa},
\]
where
\[
\Delta(u,v)=\prod_{i=1}^{l(u)}\prod_{j=1}^{l(v)} (u_i-v_j)
\]
for $u=(u_1,\dots,u_{l(u)})$ and $v=(v_1,\dots,v_{l(v)})$.
In the case of $\gsl_2$
the phase function coincides with the integrand of the Selberg integral
after identifying $\gamma=1/\kappa$, $\alpha=1-\la_1/\kappa$ and
$\beta=1-\mu_1/\kappa$.

In \cite{TV03} Tarasov and Varchenko
dealt with the $\A{2}$ case of \eqref{MV},
obtaining a closed form evaluation for
$\la=\la_2\La_2$ and $\mu=\mu_1\La_1+\mu_2\La_2$.
In the present paper we utilise
the theory of Macdonald polynomials to extend this to
$\A{n}$, and one of our main results is an explicit evaluation of
\eqref{MV} for $\la=\la_n\La_n$ and $\mu=\sum_i \mu_i\La_i$.
If we write $\kappa=1/\gamma$, $\la_i=(1-\alpha_i)/\gamma$ 
(so that $\alpha_1=\dots=\alpha_{n-1}=1$) and
$\mu_i=(1-\beta_i)/\gamma$,
and let $\Delta=C^{k_1,\dots,k_n}_{\gamma}[0,1]$ be the integration domain 
defined in \eqref{chain} of Section~\ref{sec4}, we may claim an 
evaluation of the \eqref{MV} for $\g=\A{n}$.
\begin{theorem}[$\A{n}$ Selberg integral]\label{thmSelbergAn}
For $n$ a positive integer let $0\leq k_1\leq k_2\leq \dots\leq k_n$ 
be integers and $k_0=k_{n+1}=0$.
Let $\alpha,\beta_1,\dots,\beta_n,\gamma\in\Complex$ such that
\begin{equation*}
\Re(\alpha)>0,~\Re(\beta_1)>0,\dots,\Re(\beta_n)>0,
\end{equation*}
\begin{equation*}
-\min\{\Re(\alpha)/(k_n-1),1/k_n\}<\Re(\gamma)<1/k_n
\end{equation*}
and
\begin{equation*}
-\Re(\beta_s)/(k_s-k_{s-1}-1)<\Re(\gamma)<
\Re(\beta_s+\cdots+\beta_r)/(r-s)
\end{equation*}
for $1\leq s\leq r\leq n$. Then
\begin{align*}
&\Int_{C^{k_1,\dots,k_n}_{\gamma}[0,1]}
\prod_{s=1}^n 
\biggl[ \Abs{\Delta\bigl(t^{(s)}\bigr)}^{2\gamma}
\prod_{i=1}^{k_s} \bigl(t_i^{(s)}\bigr)^{\alpha_s-1}
\bigl(1-t_i^{(s)}\bigr)^{\beta_s-1}\biggr] 
\prod_{s=1}^{n-1}
\Abs{\Delta\bigl(t^{(s)},t^{(s+1)}\bigr)}^{-\gamma} 
\; \dup t  \\ 
&=\prod_{1\leq s\leq r\leq n} \prod_{i=1}^{k_s-k_{s-1}}
\frac{\Gamma(\beta_s+\cdots+\beta_r+(i+s-r-1)\gamma)}
{\Gamma(\alpha_r
+\beta_s+\cdots+\beta_r+(i+s-r+k_r-k_{r+1}-2)\gamma)} \\
&\quad\times
\prod_{s=1}^n \prod_{i=1}^{k_s}
\frac{\Gamma(\alpha_s+(i-k_{s+1}-1)\gamma)
\Gamma(i\gamma)}{\Gamma(\gamma)},
\end{align*}
where $\alpha_1=\cdots=\alpha_{n-1}=1$, $\alpha_n=\alpha$
and $\dup t=\dup t^{(1)}\cdots\dup t^{(n)}$.
\end{theorem}

\noindent\textbf{Remarks.}

\begin{enumerate}
\item
Whenever $A/0$ occurs in the conditions 
on $\alpha$, $\beta_1,\dots,\beta_n$ and $\gamma$
this is to be interpreted as $\pm\infty$ with the
sign that of $A$.
This ensures the conditions are correct provided
$k_s>k_{s-1}$ for all $1\leq s\leq n$. 
Only minor modifications are required if
$k_s=k_{s-1}$ for some $s$.
We also note that the condition $\Re(\gamma)<1/k_n$ 
comes from $\Re(\gamma)<\min\{1/k_s: 2\leq s\leq n\}$
and does not apply when $n=1$.
\item

For $k_1=\dots=k_{n-1}=0$ and $(k_n,\beta_n,t^{(n)})\to (k,\beta,t)$ the 
$\A{n}$ Selberg integral simplifies to
\begin{multline*}
\qquad\quad
\Int_{C^{0,\dots,0,k}_{\gamma}[0,1]}\abs{\Delta(t)}^{2\gamma}
\prod_{i=1}^k
t_i^{\alpha-1}(1-t_i)^{\beta-1} \dup t \\
=\prod_{i=1}^k
\frac{\Gamma(\alpha+(i-1)\gamma)
\Gamma(\beta+(i-1)\gamma)\Gamma(i\gamma)}
{\Gamma(\alpha+\beta+(i+k-2)\gamma)\Gamma(\gamma)} .
\end{multline*}
Since (see \eqref{chain})
\begin{equation*}
C^{0,\dots,0,k}_{\gamma}[0,1]=\{t\in \R^k\,|\,
0\leq t_k\leq t_{k-1}\leq \dots\leq t_1\leq 1\}
\end{equation*}
this is equivalent to the Selberg integral \eqref{Selberg}.
Indeed, by the symmetry of the integrand 
we may replace $C^{0,\dots,0,k}_{\gamma}[0,1]$
by $[0,1]^k$ provided the right-hand side is multiplied by $k!$. 
Absorbing this factor in the ratio of gamma functions 
yields \eqref{Selberg}.
More generally the integration domain 
$C^{k_1,\dots,k_n}_{\gamma}[0,1]$ is such that $\Delta(t^{(s)})\geq 0$,
and the absolute value sign in $\abs{\Delta(t^{(s)})}^{2\gamma}$
(but not in $\abs{\Delta(t^{(s)},t^{(s+1)})}^{-\gamma}$) may be omitted.

\item Denoting the $\A{n}$ Selberg integral by
\begin{equation*}
I_{k_1,\dots,k_n}^{\A{n}}(\alpha;\beta_1,\dots,\beta_n;\gamma),
\end{equation*}
it readily follows that 
\begin{equation*}
I_{\underbrace{\scriptstyle 0,\dots,0,}_{n-m} l_1,\dots,l_m}^{\A{n}}
(\alpha;\beta_1,\dots,\beta_n;\gamma)
=I_{l_1,\dots,l_m}^{\A{m}}(\alpha;\beta_{n-m+1},\dots,\beta_n;\gamma).
\end{equation*}
In particular we have
\begin{equation*}
C^{\overbrace{\scriptstyle 0,\dots,0}^{n-m},l_1,\dots,l_m}_{\gamma}[0,1]
=C^{l_1,\dots,l_m}_{\gamma}[0,1].
\end{equation*}
The case $k_1=\dots=k_{n-1}=0$, $k_n=k$ discussed in (2) is of course a
special case of this more general reduction formula.

\item
By an appropriate change of integration variables (see Section~\ref{sec7}
for details) it follows that
\begin{multline}\label{keen}
\qquad\qquad I_{1,1,k_3,\dots,k_n}^{\A{n}}
(\alpha;\beta_1,\beta_2,\dots,\beta_n;\gamma) \\
=I_{1,k_3,\dots,k_n}^{\A{n-1}} 
(\alpha;\beta_1+\beta_2-\gamma,\beta_3,\dots,\beta_n;\gamma) \,
\frac{\Gamma(1-\gamma)\Gamma(\beta_1)}{\Gamma(\beta_1-\gamma+1)}.
\end{multline}
By iteration all but the last $k_s$ equal to $1$ may thus be eliminated.

\item
Upon taking $n=2$ and $(k_1,k_2)\to (k_2,k_1)$ and 
$(\beta_1,\beta_2)\to(\beta_2,\beta_1)$ we obtain the $\gsl_3$ Selberg
integral of Tarasov and Varchenko \cite[Theorem 3.3]{TV03}, 
see also \cite{Varchenko04}.

\item
If we denoting the set of positive roots of $\A{n}$ by $\Phi_{+}$
(i.e., $\Phi_{+}=\{\ba_s+\cdots+\ba_r:~1\leq s\leq r\leq n\}$)
then the product over ${1\leq s\leq r\leq n}$
on the right-hand side of the $\A{n}$ Selberg integral 
corresponds to the following product over $\Phi_{+}$:
\[
\prod_{1\leq s\leq r\leq n} g(\beta_s+\cdots+\beta_r)
=\prod_{\ba\in\Phi_{+}} g(\bil{\Lambda}{\ba}),
\]
where $\Lambda=\beta_1\Lambda_1+\cdots+\beta_n\Lambda_n$.

\end{enumerate}

By replacing $(\beta_s,t^{(s)})\to(\zeta\beta_s,t^{(s)}/\zeta)$ with 
$\zeta\in\R$ and then letting $\zeta$ tend to infinity we obtain the 
following exponential form of Theorem~\ref{thmSelbergAn}, with
$C^{k_1,\dots,k_n}_{\gamma}[0,\infty]$ the domain defined in \eqref{chain} 
of Section~\ref{sec4}.
\begin{corollary}[First $\A{n}$ exponential Selberg integral]\label{thmESelbergAn}
For $n$ a positive integer let $0\leq k_1\leq k_2\leq \dots\leq k_n$ 
be integers and $k_0=k_{n+1}=0$.
Let $\alpha,\beta_1,\dots,\beta_n,\gamma\in\Complex$ such that
\begin{equation*}
\Re(\alpha)>0,~\Re(\beta_1)>0,\dots,\Re(\beta_n)>0
\end{equation*}
and
\begin{equation*}
-\min\{\Re(\alpha)/(k_n-1),1/k_n\}<\Re(\gamma)<1/k_n.
\end{equation*}
Then
\begin{align*}
&\Int_{C^{k_1,\dots,k_n}_{\gamma}[0,\infty]}
\prod_{s=1}^n \biggl[ 
\Abs{\Delta\bigl(t^{(s)}\bigr)}^{2\gamma} 
\prod_{i=1}^{k_s}
\bigl(t_i^{(s)}\bigr)^{\alpha_s-1}
\eup^{-\beta_s t_i^{(s)}} \biggr]
\prod_{s=1}^{n-1}
\Abs{\Delta\bigl(t^{(s)},t^{(s+1)}\bigr)}^{-\gamma} \; \dup t \\
&\qquad =\prod_{1\leq s\leq r\leq n} 
(\beta_s+\cdots+\beta_r)^{-(\alpha_r+(k_r-k_{r+1}-1)\gamma)(k_s-k_{s-1})} \\
&\qquad\quad\times
\prod_{s=1}^n \prod_{i=1}^{k_s}
\frac{\Gamma(\alpha_s+(i-k_{s+1}-1)\gamma)\Gamma(i\gamma)}{\Gamma(\gamma)},
\end{align*}
\begin{comment}
\begin{align*}
&\Int_{C^{k_1,\dots,k_n}_{\gamma}[0,\infty]}
\prod_{i=1}^{k_n}
\bigl(t_i^{(n)}\bigr)^{\alpha-1}
\prod_{s=1}^n \biggl[\, \eup^{-\beta_s\sum_{i=1}^{k_s}t_i^{(s)}} 
\Abs{\Delta\bigl(t^{(s)}\bigr)}^{2\gamma} \,\biggr] \\
&\qquad\qquad\quad\times
\prod_{s=1}^{n-1}
\Abs{\Delta\bigl(t^{(s)},t^{(s+1)}\bigr)}^{-\gamma} 
\; \dup t^{(1)}\cdots\dup t^{(n)} \\
&\qquad =\prod_{1\leq s\leq r\leq n} 
(\beta_s+\cdots+\beta_r)^{-(\alpha_r+(k_r-k_{r+1}-1)\gamma)(k_s-k_{s-1})} \\
&\qquad\quad\times
\prod_{s=1}^n \prod_{i=1}^{k_s}
\frac{\Gamma(\alpha_s+(i-k_{s+1}-1)\gamma)\Gamma(i\gamma)}{\Gamma(\gamma)},
\end{align*}
\end{comment}
where
$\alpha_1=\cdots=\alpha_{n-1}=1$ and $\alpha_n=\alpha$.
\end{corollary}
Replacing $t_i\to 1-t_i$ for all $1\leq i\leq k$ in the Selberg 
integral \eqref{Selberg} leads to an interchange of $\alpha$ and 
$\beta$. Consequently the classical Selberg integral has just a 
single exponential form. This $(\alpha,\beta)$-symmetry is 
no longer present for $n>1$, and replacing $t_i^{(s)}\to 1-t_i^{(s)}$
for all $1\leq i\leq k_s$ and $1\leq s\leq n$ followed by
$(\alpha,t^{(s)})\to(\zeta,t^{(s)}/\zeta)$ (with $\zeta\in\R$)
and then letting $\zeta$ tend to infinity, results in a second
exponential form of Theorem~\ref{thmSelbergAn}.
Below $\bar{C}^{k_1,\dots,k_n}_{\gamma}[0,\infty]$ is the
integration domain defined in \eqref{chain2} of Section~\ref{sec5}.
\begin{corollary}[Second $\A{n}$ exponential Selberg integral]
For $n$ a positive integer let $0\leq k_1\leq k_2\leq \dots\leq k_n$
be integers and $k_0=k_{n+1}=0$.
Let $\beta_1,\dots,\beta_n,\gamma\in\Complex$ such that
\begin{equation*}
\Re(\beta_1)>0,\dots,\Re(\beta_n)>0,
\qquad -1/k_n<\Re(\gamma)<1/k_n
\end{equation*}
and
\begin{equation*}
-\Re(\beta_s)/(k_s-k_{s-1}-1)<\Re(\gamma)<
\Re(\beta_s+\cdots+\beta_r)/(r-s)
\end{equation*}
for $1\leq s\leq r\leq n$. Then
\begin{align*}
&\Int_{\bar{C}^{k_1,\dots,k_n}_{\gamma}[0,\infty]}
\!\!
\eup^{-\sum_{i=1}^{k_n} t_i^{(n)}}
\prod_{s=1}^n \biggl[
\Abs{\Delta\bigl(t^{(s)}\bigr)}^{2\gamma}
\prod_{i=1}^{k_s}\bigl(t_i^{(s)}\bigr)^{\beta_s-1}\biggr]
\prod_{s=1}^{n-1}
\Abs{\Delta\bigl(t^{(s)},t^{(s+1)}\bigr)}^{-\gamma}
\dup t \\
&=\prod_{1\leq s\leq r\leq n} \prod_{i=1}^{k_s-k_{s-1}}
\Gamma(\beta_s+\cdots+\beta_r+(i+s-r-1)\gamma)\\
&\quad\times\prod_{1\leq s\leq r\leq n-1} \prod_{i=1}^{k_s-k_{s-1}}
\frac{1}{\Gamma(1+\beta_s+\cdots+\beta_r+(i+s-r+k_r-k_{r+1}-2)\gamma)} \\
&\quad\times
\prod_{s=1}^n \prod_{i=1}^{k_s}
\frac{\Gamma(i\gamma)}{\Gamma(\gamma)}\:
\prod_{s=1}^{n-1} \prod_{i=1}^{k_s}
\Gamma(1+(i-k_{s+1}-1)\gamma).
\end{align*}
\end{corollary}

\medskip

\subsection{Outline}
In Section~\ref{sec2} we review some standard facts about Macdonald 
polynomials needed to prove an identity for the $q,t$-analogues of the
Littlewood--Richard\-son coefficients (Theorem~\ref{thm1}). 
In Section~\ref{sec3} we apply Theorem~\ref{thm1} to establish a new $\A{n}$
$q$-binomial theorem for Macdonald polynomials (Theorem~\ref{thm3}).
In Section~\ref{sec4} we first utilise the $q=1$ case of this 
theorem to prove the exponential $\A{n}$ Selberg integral of 
Corollary~\ref{thmESelbergAn}.
Then, in Section~\ref{sec5}, we exploit the full $q$-binomial theorem 
to obtain a multidimensional $q$-integral which yields 
Theorem~\ref{thmSelbergAn} in the $q\to 1$ limit.
In Section~\ref{sec6} we generalise the $\A{n}$ Selberg integral
by including a Jack polynomial in the integrand (Theorem~\ref{thmJack}).
Finally, in Section~\ref{sec7}, we give the full details of
two special cases of Theorem~\ref{thmSelbergAn}, corresponding to
$(k_1,\dots,k_{n-1},k_n)=(1,\dots,1,k)$ and $\gamma=0$ respectively.

\section{Macdonald polynomials}\label{sec2}
Our main tool in the proof of Theorem~\ref{thmSelbergAn} is 
the theory of symmetric functions, and in Sections~\ref{secpr} and
\ref{secMD} we review some well-known facts from the theory. 
For a more comprehensive
introduction we refer the reader to~\cite{Lascoux03,Macdonald95,Stanley99}.

\subsection{Preliminaries}\label{secpr}
Let $\la=(\la_1,\la_2,\dots)$ be a partition,
i.e., $\la_1\geq \la_2\geq \dots$ with finitely many $\la_i$
unequal to zero.
The length and weight of $\la$, denoted by
$l(\la)$ and $\abs{\la}$, are the number and sum
of the non-zero $\la_i$ respectively.
As usual we identify two partitions that differ only in their
string of zeros, so that $(6,3,3,1,0,0)$ and $(6,3,1,1)$ represent the
same partition.
When $\abs{\la}=N$ we say that $\la$ is a partition of $N$,
and the unique partition of zero is denoted by $0$.
The multiplicity of the part $i$ in the partition $\la$ is denoted
by $m_i=m_i(\la)$, and occasionally we will write
$\la=(1^{m_1} 2^{m_2} \dots)$.

We identify a partition with its Ferrers graph,
defined by the set of points in $(i,j)\in \Z^2$ such that
$1\leq j\leq \la_i$, and further make the usual 
identification between Ferrers graphs and (Young) diagrams
by replacing points by squares.

The conjugate $\la'$ of $\la$ is the partition obtained by
reflecting the diagram of $\la$ in the main diagonal,
so that, in particular, $m_i(\la)=\la_i'-\la_{i+1}'$.
The statistic $n(\la)$ is given by
\begin{equation*}
n(\la)=\sum_{i\geq 1} (i-1)\la_i=
\sum_{i\geq 1}\binom{\la_i'}{2}.
\end{equation*}

The dominance partial order on the set of partitions of $N$ is
defined by $\la\geq \mu$ if
$\la_1+\cdots+\la_i\geq \mu_1+\cdots+\mu_i$ for all $i\geq 1$.
If $\la\geq \mu$ and $\la\neq\mu$ then $\la>\mu$.

If $\la$ and $\mu$ are partitions then $\mu\subseteq\la$
if (the diagram of) $\mu$ is contained in (the diagram of)
$\la$, i.e., $\mu_i\leq\la_i$ for all $i\geq 1$.
If $\mu\subseteq\la$ then the skew-diagram $\la-\mu$
denotes the set-theoretic difference between $\la$ and $\mu$,
i.e., those squares of $\la$ not contained in $\mu$.

Let $s=(i,j)$ be a square in the diagram of $\la$. Then
$a(s)$, $a'(s)$, $l(s)$ and $l'(s)$ are the arm-length, arm-colength,
leg-length and leg-colength of $s$, defined by
\begin{subequations}\label{aapllp}
\begin{align}
a(s)&=\la_i-j,  & a'(s)&=j-1 \\
l(s)&=\la'_j-i,  & l'(s)&=i-1.
\end{align}
\end{subequations}
This may be used to define the generalised hook-length polynomials
\cite[Equation (VI.8.1)]{Macdonald95}
\begin{subequations}\label{cdef}
\begin{align}
c_{\la}(q,t)&=\prod_{s\in\la}\bigl(1-q^{a(s)}t^{l(s)+1}\bigr), \\
c'_{\la}(q,t)&=\prod_{s\in\la}\bigl(1-q^{a(s)+1}t^{l(s)}\bigr),
\end{align}
\end{subequations}
where the products are over all squares of $\la$. We further set
\begin{equation}\label{bdef}
b_{\la}(q,t)=\frac{c_{\la}(q,t)}{c'_{\la}(q,t)}.
\end{equation}

For $N$ a nonnegative integer the $q$-shifted factorial 
$(b;q)_N$ is defined as $(b;q)_0=1$ and
\begin{equation}\label{qfac}
(b;q)_N=(1-b)(1-bq)\cdots(1-bq^{N-1}).
\end{equation}
We also need the $q$-shifted factorial for negative 
(integer) values of $N$. This may be obtained from the above by
\begin{equation*}
(b;q)_{-N}=\frac{1}{(bq^{-N};q)_N}.
\end{equation*}
This implies in particular that $1/(q;q)_{-N}=0$ for positive $N$.

The definition \eqref{qfac} may be extended to partitions by
\begin{equation*}
(b;q,t)_{\la}=\prod_{s\in\la}\bigl(1-b q^{a'(s)}t^{-l'(s)}\bigr) \\
=\prod_{i=1}^{l(\la)}(bt^{1-i};q)_{\la_i}.
\end{equation*}
With this notation the polynomials \eqref{cdef} can be expressed as
\cite[Proposition 3.2]{Kaneko96}
\begin{subequations}\label{ccp}
\begin{align}\label{c}
c_{\la}(q,t)&=(t^n;q,t)_{\la}
\prod_{1\leq i<j\leq n}\frac{(t^{j-i};q)_{\la_i-\la_j}}
{(t^{j-i+1};q)_{\la_i-\la_j}}, \\
\label{cp}
c'_{\la}(q,t)&=(qt^{n-1};q,t)_{\la}
\prod_{1\leq i<j\leq n}\frac{(qt^{j-i-1};q)_{\la_i-\la_j}}
{(qt^{j-i};q)_{\la_i-\la_j}},
\end{align}
\end{subequations}
where $n$ is any integer such that $n\geq l(\la)$.

Finally we introduce the usual condensed notation for $q$-shifted factorials
as
\begin{equation*}
(a_1,\dots,a_k;q)_N=(a_1;q)_N\cdots (a_k;q)_N
\end{equation*}
and
\begin{equation*}
(a_1,\dots,a_k;q,t)_{\la}=(a_1;q,t)_{\la}\cdots (a_k;q,t)_{\la}.
\end{equation*}

\subsection{Macdonald polynomials}\label{secMD}
Let $\Symm_n$ denote the symmetric group,
and $\Lambda_n=\Z[x_1,\dots,x_n]^{\Symm_n}$
the ring of symmetric polynomials in $n$ independent variables.

For $x=(x_1,\dots,x_n)$ and $\la=(\la_1,\dots,\la_n)$ a partition 
of at most $n$ parts the monomial symmetric function $m_{\la}$ 
is defined as
\begin{equation*}
m_{\la}(x)=\sum x^{\alpha}.
\end{equation*}
Here the sum is over all distinct permutations $\alpha$ of
$\la$, and $x^{\alpha}=x_1^{\alpha_1}\cdots x_n^{\alpha_n}$.
For $l(\la)>n$ we set $m_{\la}(x)=0$.
The monomial symmetric functions $m_{\la}$ for $l(\la)\leq n$
form a $\Z$-basis of $\Lambda_n$.

For $r$ a nonnegative integer the power sums $p_r$ are given by
$p_0=1$ and $p_r=m_{(r)}$ for $r>1$. Hence
\begin{equation*}
p_r(x)=\sum_{i=1}^n x_i^r.
\end{equation*}
More generally the power-sum products are defined as
$p_{\la}(x)=p_{\la_1}(x)\cdots p_{\la_n}(x)$.

Following Macdonald we define the scalar product
$\langle \cdot,\cdot \rangle_{q,t}$ by
\begin{equation*}
\langle p_{\la},p_{\mu}\rangle_{q,t}=
\delta_{\la\mu} z_{\la} \prod_{i=1}^n
\frac{1-q^{\la_i}}{1-t^{\la_i}},
\end{equation*}
with $z_{\la}=\prod_{i\geq 1} m_i! \: i^{m_i}$ and
$m_i=m_i(\la)$. If we denote the ring of symmetric functions 
in $n$ variables over the field $\F=\Q(q,t)$ of rational functions 
in $q$ and $t$ by $\Lambda_{n,\F}$, then
the Macdonald polynomial $P_{\la}(x;q,t)$
is the unique symmetric polynomial in $\Lambda_{n,\F}$ such that
\cite[Equation (VI.4.7)]{Macdonald95}:
\begin{equation}\label{Pm}
P_{\la}(x;q,t)=m_{\la}(x)+\sum_{\mu<\la}
u_{\la\mu}(q,t) m_{\mu}(x)
\end{equation}
and
\begin{equation*}
\langle P_{\la},P_{\mu} \rangle_{q,t}
=0\quad \text{if$\quad\la\neq\mu$.}
\end{equation*}
The Macdonald polynomials $P_{\la}(x;q,t)$ with $l(\la)\leq n$ 
form an $\F$-basis of $\Lambda_{n,\F}$. If $l(\la)>n$ then
$P_{\la}(x;q,t)=0$.
{}From \eqref{Pm} it follows that $P_{\la}(x;q,t)$ for
$l(\la)\leq n$ is homogeneous of degree $\abs{\la}$:
\begin{equation}\label{hom}
P_{\la}(zx;q,t)=z^{\abs{\la}}P_{\la}(x;q,t)
\end{equation}
with $z$ a scalar.

For $f\in\Lambda_{n,\F}$ and $\la$ a partition such that $l(\la)\leq n$
the evaluation homomorphism 
$u_{\la}^{(n)}:\Lambda_{n,\F}\to \F$ is defined as
\begin{equation}\label{eval}
u_{\la}^{(n)}(f)=f(q^{\la_1}t^{n-1},q^{\la_2}t^{n-2},\dots,q^{\la_n}t^0).
\end{equation}
We extend this to 
$f\in\F(x_1,\dots,x_n)^{\Symm_n}$ for those $f$ for which the
right-hand side of \eqref{eval} is well-defined.
The principal specialisation formula for Macdonald polynomials 
corresponds to
\cite[Example VI.6.5]{Macdonald95}
\begin{equation}\label{PS}
u_0^{(n)}(P_{\la})
=t^{n(\la)} \prod_{s\in\la}
\frac{1-q^{a'(s)}t^{n-l'(s)}}{1-q^{a(s)}t^{l(s)+1}}
=t^{n(\la)}\frac{(t^n;q,t)_{\la}}{c_{\la}(q,t)}.
\end{equation}
For more general evaluations we have the symmetry
\cite[Equation (VI.6.6)]{Macdonald95}
\begin{equation}\label{symm}
u_{\la}^{(n)}(P_\mu)u_0^{(n)}(P_{\la})=
u_{\mu}^{(n)}(P_{\la})u_0^{(n)}(P_{\mu})
\end{equation}
for $l(\la),l(\mu)\leq n$. 
It will be convenient to also define $u_{\la;z}^{(n)}$ as
\begin{equation}\label{evalz}
u_{\la;z}^{(n)}(f)=
f(zq^{\la_1}t^{n-1},zq^{\la_2}t^{n-2},\dots,zq^{\la_n}t^0).
\end{equation}
For homogeneous functions of degree $d$ we of course have
\begin{equation}\label{deg}
u_{\la;z}^{(n)}(f)=z^d \, u_{\la}^{(n)}(f).
\end{equation}

Thanks to the stability $P_{\la}(x_1,\dots,x_n;q,t)=
P_{\la}(x_1,\dots,x_n,0;q,t)$ for $l(\la)\leq n$,
we may extend the $P_{\la}$ to an infinite alphabet, and
in the remainder of this section
we assume that $x$ (and $y$) contain countably many variables.
By abuse of terminology we will still refer to $P_{\la}(x;q,t)$
as a Macdonald polynomial, instead of a Macdonald function.
Then the Cauchy identity for Macdonald polynomials
is given by \cite[Equation (VI.4.13)]{Macdonald95}
\begin{equation}\label{Cauchy}
\sum_{\la} b_{\la}(q,t)P_{\la}(x;q,t)P_{\la}(y;q,t)
=\prod_{i,j\geq 1} \frac{(tx_i y_j;q)_{\infty}}{(x_i y_j;q)_{\infty}},
\end{equation}
with $b_{\la}(q,t)$ defined in \eqref{bdef}.

The $q,t$-Littlewood--Richardson coefficients are defined as
\begin{equation}\label{qtLR}
P_{\mu}(x;q,t)P_{\nu}(x;q,t)=
\sum_{\la}f_{\mu\nu}^{\la}(q,t) P_{\la}(x;q,t)
\end{equation}
and trivially satisfy
\begin{equation*}
f_{\mu\nu}^{\la}(q,t)=f_{\nu\mu}^{\la}(q,t)
\end{equation*}
and
\begin{equation}\label{SA}
f_{\mu\nu}^{\la}(q,t)=0 \text{ unless }\abs{\la}=\abs{\mu}+\abs{\nu}.
\end{equation}
It can also be shown that 
\cite[Equation (VI.7.7)]{Macdonald95}
\begin{equation}\label{include}
f_{\mu\nu}^{\la}(q,t)=0 \text{ unless } \mu,\nu\subseteq\la.
\end{equation}

The $q,t$-Littlewood--Richardson coefficients may be used to define
the skew Macdonald polynomials
\begin{equation}\label{skewdef}
P_{\la/\mu}(x;q,t)=\sum_{\nu}f_{\mu\nu}^{\la}(q,t)P_{\nu}(x;q,t).
\end{equation}
By \eqref{include}, $P_{\la/\mu}(x;q,t)=0$
unless $\mu\subseteq\la$ (in which case it is a homogeneous
of degree $\abs{\la}-\abs{\mu}$.
Equivalent to \eqref{skewdef} is
\begin{equation}\label{skewdef2}
P_{\la}(x,y;q,t)=\sum_{\mu} P_{\la/\mu}(x;q,t)P_{\mu}(y;q,t).
\end{equation}

Finally we need the Kaneko--Macdonald definition
of basic hypergeometric series with Macdonald polynomial
argument \cite{Kaneko96,Macdonald}
\begin{equation}\label{Phirs}
{_{r+1}\Phi_r}\biggl[\genfrac{}{}{0pt}{}
{a_1,\dots,a_{r+1}}{b_1,\dots,b_r};q,t;x\biggr]
=\sum_{\lambda} 
t^{n(\la)}\frac{P_{\lambda}(x;q,t)}{c'_{\lambda}(q,t)}\,
\frac{(a_1,\dots,a_{r+1};q,t)_{\lambda}}
{(b_1,\dots,b_r;q,t)_{\lambda}}.
\end{equation}
When $x=(z)$ this reduces to the classical $_{r+1}\phi_r$
basic hypergeometric series \cite{GR04}:
\begin{align*}
{_{r+1}}\Phi_r\biggl[\genfrac{}{}{0pt}{}
{a_1,\dots,a_{r+1}}{b_1,\dots,b_r};q,t;(z)\biggr]
&=\sum_{k=0}^{\infty}
\frac{(a_1,\dots,a_{r+1};q)_k}
{(q,b_1,\dots,b_r;q)_k} \, z^k  \\
&={_{r+1}}\phi_r\biggl[\genfrac{}{}{0pt}{}
{a_1,\dots,a_{r+1}}{b_1,\dots,b_r};q,z\biggr].
\end{align*}
The main result needed for $_{r+1}\Phi_r$ series is
the $q$-binomial theorem
\cite[Theorem 3.5]{Kaneko96}, \cite[Equation (2.2)]{Macdonald}
(see also \cite[Theorem 3]{Kaneko96b} and \cite[Lemma 3.1]{MN98})
\begin{equation}\label{Phi10sum}
{_1\Phi_0}\biggl[\genfrac{}{}{0pt}{}
{a}{\text{--}}\,;q,t;x\biggr]
=\prod_{i\geq 1} \frac{(ax_i;q)_{\infty}}{(x_i;q)_{\infty}}.
\end{equation}

\noindent\textbf{Remark.}
In this paper we mostly view results such as
\eqref{Cauchy} and \eqref{Phi10sum} as formal identities.
Later, when transforming formal power series to
integrals, issues of convergence do become important.
It is however not difficult to give necessary convergence
conditions for each of the identities in this paper.
For example, in \eqref{Phi10sum}, we may add $x=(x_1,\dots,x_n)$,
$\abs{q}<1$ and $\max\{\abs{x_1},\dots,\abs{x_n}\}<1$,
and view the $_1\Phi_0$ as a genuine hypergeometric function.

\subsection{An identity for $q,t$-Littlewood--Richardson coefficients}
\label{secqtLR}
A crucial role in our proof of the $\A{n}$ Selberg integral of
Theorem~\ref{thmSelbergAn} is the following identity for the 
$q,t$-Littlewood--Richardson coefficients.
\begin{theorem}\label{thm1}
Given two integers $0\leq m\leq n$, let $\la$ and $\mu$ be 
partitions such that $l(\la)\leq m$ and $l(\mu)\leq n$.
Then
\begin{multline}\label{mn}
\sum_{\nu,\omega}
 t^{n(\nu)-\abs{\omega}}
f_{\omega\nu}^{\la}(q,t)u_0^{(n-m)}(P_{\mu/\omega})\,
\frac{(qt^{m-n-1};q,t)_{\nu}}{c'_{\nu}(q,t)} \\
=t^{n(\la)-m\abs{\mu}}
u_0^{(n)}(P_{\mu})\,
\frac{(qt^{m-1};q,t)_{\la}}{c'_{\la}(q,t)}
\prod_{i=1}^m
\prod_{j=1}^n
\frac{(qt^{j-i+m-n-1};q)_{\la_i-\mu_j}}
{(qt^{j-i+m-n};q)_{\la_i-\mu_j}}.
\end{multline}
\end{theorem}
Since $f_{\omega\nu}^{\la}(q,t)=0$ if $\omega\not\subseteq\la$ and
$P_{\mu/\omega}=0$ if $\omega\not\subseteq\mu$ we may add the
restrictions $\omega\subseteq\la$ and $\omega\subseteq\mu$ to
the sum over $\omega$.
We will in fact show that the summand on the left vanishes unless
\begin{equation}\label{van}
\la_i\geq \mu_{i+n-m} \quad\text{for~~$1\leq i\leq m$.}
\end{equation}
In other words, if $\mu^{\ast}$ is the partition formed by the
last $m$ parts of $\mu$ (i.e., $\mu^{\ast}=(\mu_{n-m+1},\dots,\mu_n)$)
then the summand vanishes unless $\mu^{\ast}\subseteq\la$.
To see this we recall from \cite[Equation (VI.7.13$'$)]{Macdonald95} that
\begin{equation*}
P_{\mu/\omega}(x_1,\dots,x_{n-m};q,t)=\sum_{T}\psi_T(q,t)x^T,
\end{equation*}
where the sum is over all semistandard Young tableaux $T$ of 
skew shape $\mu-\omega$ over the alphabet $\{1,\dots,n-m\}$;
$x^T$ is the monomial defined by $T$ and $\psi_T\in \F$.
For the shape $\mu-\omega$ to have an admissible filling it
must have at most $n-m$ boxes in each of its columns.
Hence $\omega_i\geq\mu_{i+n-m}$ for $1\leq i\leq m$.
Since we already established that the summand vanishes unless
$\omega\subseteq\la$, a necessary condition for nonvanishing of
the summand is thus given by \eqref{van}.
Since $1/(q;q)_{-N}=0$ for $N$ a positive integer,
it is easily seen that also the double product
on the right-hand side of \eqref{mn} vanishes unless
\eqref{van} holds.

Theorem~\ref{thm1} for arbitary $0\leq m\leq n$ corresponds to the $u=0$ case 
of a more general result established in \cite[Theorem 4.1]{Warnaar05}.
For $m=n$, so that 
\[
P_{\mu/\omega}(x_1,\dots,x_{n-m};q,t)=\delta_{\mu\omega},
\]
the theorem simplifies to \cite[Proposition 3.2]{Warnaar}.
A proof Theorem~\ref{thm1} is included below for the sake of completeness.

\begin{proof}[Proof of Theorem~\ref{thm1}]
Let $x=(x_1,\dots,x_m)$ and $y=(y_1,\dots,y_n)$ so that
the Cauchy identity \eqref{Cauchy} becomes
\begin{equation*}
\sum_{\eta} b_{\eta}(q,t)P_{\eta}(x;q,t)P_{\eta}(y;q,t)
=\prod_{i=1}^m \prod_{j=1}^n 
\frac{(tx_i y_j;q)_{\infty}}{(x_i y_j;q)_{\infty}}.
\end{equation*}
Next we apply the homomorphisms 
$u_{\la;z}^{(m)}$ (acting on $x$) and $u_{\mu}^{(n)}$ (acting on $y$),
and use the homogeneity \eqref{hom} of the Macdonald polynomials. Hence
\begin{multline}\label{eta}
\sum_{\eta} 
z^{\abs{\eta}} b_{\eta}(q,t)u_{\la}^{(m)}(P_{\eta})u_{\mu}^{(n)}(P_{\eta}) \\
=\prod_{i=1}^m \frac{(zt^{n+m-i};q)_{\infty}}{(zt^{m-i};q)_{\infty}}
\prod_{i=1}^m \prod_{j=1}^n 
\frac{(zt^{n+m-i-j};q)_{\la_i+\mu_j}}{(zt^{n+m-i-j+1};q)_{\la_i+\mu_j}}.
\end{multline}
The summand on the left vanishes unless $l(\eta)\leq\min\{n,m\}$.
Assuming such $\eta$ we may twice use
the symmetry \eqref{symm} to rewrite the left-hand side as
\begin{equation*}
\text{LHS}\eqref{eta}=
\sum_{\eta} z^{\abs{\eta}} b_{\eta}(q,t)\,
\frac{u_{\eta}^{(m)}(P_{\la})u_{\eta}^{(n)}(P_{\mu}) 
u_0^{(m)}(P_{\eta})u_0^{(n)}(P_{\eta})}
{u_0^{(m)}(P_{\la})u_0^{(n)}(P_{\mu})}.
\end{equation*}
In the remainder we assume that $n\geq m$ and apply \eqref{skewdef2} as well as
\eqref{hom} to get
\begin{align*}
u_{\eta}^{(n)}(P_{\mu})&=
P_{\mu}(q^{\eta_1}t^{n-1},\dots,q^{\eta_m}t^{n-m},t^{n-m-1},\dots,t,1;q,t) \\[1mm]
&=\sum_{\omega}P_{\omega}(q^{\eta_1}t^{n-1},\dots,q^{\eta_m}t^{n-m};q,t)
u_0^{(n-m)}(P_{\mu/\omega}) \\
&=\sum_{\omega}t^{(n-m)\abs{\omega}}u_{\eta}^{(m)}(P_{\omega})u_0^{(n-m)}(P_{\mu/\omega}).
\end{align*}
Thus
\begin{multline*}
\text{LHS}\eqref{eta}=\sum_{\eta,\omega} 
z^{\abs{\eta}} t^{(n-m)\abs{\omega}} b_{\eta}(q,t) \\
\times
\frac{u_0^{(n-m)}(P_{\mu/\omega})u_{\eta}^{(m)}(P_{\la}) u_{\eta}^{(m)}(P_{\omega})
u_0^{(m)}(P_{\eta})u_0^{(n)}(P_{\eta})}{u_0^{(m)}(P_{\la})u_0^{(n)}(P_{\mu})}.
\end{multline*}
Next we use that
\begin{align*}
u_{\eta}^{(m)}(P_{\la}) u_{\eta}^{(m)}(P_{\omega})&=
u_{\eta}^{(m)}(P_{\la} \, P_{\omega}) \\[2mm]
&=u_{\eta}^{(m)}\Bigl(\, \sum_{\nu} f_{\omega\la}^{\nu}P_{\nu}\Bigr) 
\qquad (\text{by \eqref{qtLR})} \\
&=\sum_{\nu} f_{\omega\la}^{\nu}\, u_{\eta}^{(m)}(P_{\nu}) 
\end{align*}
to rewrite this as
\begin{multline*}
\text{LHS}\eqref{eta}=\sum_{\eta,\omega,\nu} 
z^{\abs{\eta}} t^{(n-m)\abs{\omega}}f_{\omega\la}^{\nu}(q,t) b_{\eta}(q,t) \\
\times
\frac{u_0^{(n-m)}(P_{\mu/\omega})u_{\eta}^{(m)}(P_{\nu}) 
u_0^{(m)}(P_{\eta})u_0^{(n)}(P_{\eta})}{u_0^{(m)}(P_{\la})u_0^{(n)}(P_{\mu})}.
\end{multline*}
By one more application of \eqref{symm} this becomes
\begin{multline*}
\text{LHS}\eqref{eta}=\sum_{\eta,\omega,\nu} 
z^{\abs{\eta}} t^{(n-m)\abs{\omega}}f_{\omega\la}^{\nu}(q,t) b_{\eta}(q,t) \\
\times
\frac{u_0^{(n-m)}(P_{\mu/\omega})
u_{\nu}^{(m)}(P_{\eta}) u_0^{(m)}(P_{\nu})
u_0^{(n)}(P_{\eta})}{u_0^{(m)}(P_{\la})u_0^{(n)}(P_{\mu})}.
\end{multline*}
The sum over $\eta$ may now be evaluated as follows:
\begin{align*}
\sum_{\eta} &
z^{\abs{\eta}} b_{\eta}(q,t) u_0^{(n)}(P_{\eta}) 
u_{\nu}^{(m)}(P_{\eta}) && \\
&=\sum_{\eta} t^{n(\eta)}
\frac{(t^n;q,t)_{\eta}}{c'_{\eta}(q,t)}\,
u_{\nu;z}^{(m)}(P_{\eta}) && \text{(by \eqref{bdef} and \eqref{PS})} \\
&=u_{\nu;z}^{(m)}\biggl(\sum_{\eta} t^{n(\eta)}
\frac{(t^n;q,t)_{\eta}}{c'_{\eta}(q,t)}
P_{\eta}(x;q,t)\biggr) && \\ 
&=u_{\nu;z}^{(m)}\biggl(
{_1}\Phi_0\biggl[\genfrac{}{}{0pt}{}{t^n}{\text{--}}\,;q,t;x\biggr]\biggr) 
&& (\text{by \eqref{Phirs}}) \\
&=u_{\nu;z}^{(m)}\biggl(\;
\prod_{i=1}^m \frac{(x_it^n;q)_{\infty}}{(x_i;q)_{\infty}}\biggr)
&& (\text{by \eqref{Phi10sum}}) \\
&=\prod_{i=1}^m \frac{(zq^{\nu_i}t^{n+m-i};q)_{\infty}}
{(zq^{\nu_i}t^{m-i};q)_{\infty}} && \\
&=\frac{(zt^{m-1};q,t)_{\nu}}{(zt^{n+m-1};q,t)_{\nu}}
\prod_{i=1}^m \frac{(zt^{n+m-i};q)_{\infty}}
{(zt^{m-i};q)_{\infty}}. &&
\end{align*}
We thus arrive at
\begin{multline*}
\text{LHS}\eqref{eta}=
\prod_{i=1}^m \frac{(zt^{n+m-i};q)_{\infty}}
{(zt^{m-i};q)_{\infty}} \\
\times\sum_{\omega,\nu} 
t^{(n-m)\abs{\omega}}f_{\omega\la}^{\nu}(q,t) \,
\frac{u_0^{(n-m)}(P_{\mu/\omega})
u_0^{(m)}(P_{\nu})}{u_0^{(m)}(P_{\la})u_0^{(n)}(P_{\mu})}
\frac{(zt^{m-1};q,t)_{\nu}}{(zt^{n+m-1};q,t)_{\nu}}.
\end{multline*}
Finally equating this with the right-hand side of \eqref{eta}
and replacing $z\to zt^{1-n-m}$ yields
\begin{multline}\label{zprop}
\sum_{\omega,\nu} 
t^{(n-m)\abs{\omega}}f_{\omega\la}^{\nu}(q,t) 
u_0^{(n-m)}(P_{\mu/\omega})
u_0^{(m)}(P_{\nu}) \,
\frac{(zt^{-n};q,t)_{\nu}}{(z;q,t)_{\nu}} \\
=u_0^{(m)}(P_{\la})u_0^{(n)}(P_{\mu})
\prod_{i=1}^m \prod_{j=1}^n 
\frac{(zt^{1-i-j};q)_{\la_i+\mu_j}}{(zt^{2-i-j};q)_{\la_i+\mu_j}}.
\end{multline}
Both sides of this identity trivially vanish if
$l(\la)>m$. Furthermore, the summand on the left
vanishes if $l(\nu)>m$.
Hence we may without loss of generality assume in
the following that $l(\la)\leq m$ and 
$l(\nu)\leq m$. (The latter of course refers to
a restriction on the summation index.)  
We may also assume that the largest part of $\nu$ is bounded since
$f_{\omega\la}^{\nu}=0$ if $\abs{\omega}+\abs{\la}\neq\abs{\nu}$
and $P_{\mu/\omega}=0$ if $\omega\not\subseteq\mu$.
In particular $\nu_1\leq \abs{\la}+\abs{\mu}$.

The above considerations imply that 
$\la,\nu\subseteq (N^m)$ for sufficiently large $N$.
Given such $N$ we can define the complements of $\la$ and $\nu$ with
respect to $(N^m)$. Denoting these partitions by $\hat{\la}$ and
$\hat{\mu}$ we have
$\hat{\la}_i=N-\la_{m+1-i}$ and
$\hat{\nu}_i=N-\nu_{m+1-i}$ for $1\leq i\leq m$.

We now replace $\la$ and $\nu$ by $\hat{\la}$ and $\hat{\nu}$
in \eqref{zprop} and then eliminate the hats using
\cite[page 263]{Warnaar}
\begin{equation*}
f_{\omega\hat{\la}}^{\hat{\nu}}(q,t)=
t^{n(\nu)-n(\la)}
f_{\omega\nu}^{\la}(q,t)\,
\frac{(qt^{m-1};q,t)_{\nu}}{(qt^{m-1};q,t)_{\la}}\,
\frac{c'_{\la}(q,t)}{c'_{\nu}(q,t)}\,
\frac{u_0^{(m)}(P_{\la})}{u_0^{(m)}(P_{\nu})},
\end{equation*}
\cite[Equation (4.1)]{BF99}
\begin{equation*}
(a;q,t)_{\hat{\la}}=
(-q/a)^{\abs{\la}}t^{(m-1)\abs{\la}-n(\la)}
q^{n(\la')-N\abs{\la}}
\frac{(a;q,t)_{(N^m)}}{(q^{1-N}t^{m-1}/a;q,t)_{\la}},
\end{equation*}
and
\begin{equation*}
u_0^{(m)}(P_{\hat{\la}})=t^{\binom{m}{2}N+(1-m)\abs{\la}}
u_0^{(m)}(P_{\la}).
\end{equation*}
This last result follows from \cite[Equation (4.3)]{BF99}
\begin{equation*}
P_{\hat{\la}}(x;q,t)=(x_1\cdots x_m)^N\, P_{\la}(x^{-1};q,t)
\end{equation*}
and the homogeneity \eqref{hom}.
As a result we end up with
\begin{multline*}
\sum_{\omega,\nu} 
t^{n(\nu)-\abs{\omega}}f_{\omega\nu}^{\la}(q,t) 
u_0^{(n-m)}(P_{\mu/\omega})\,
\frac{(qt^{m-1},q^{1-N}t^{m-1}/z;q,t)_{\nu}}{c'_{\nu}(q,t)\,
(q^{1-N}t^{n+m-1}/z;q,t)_{\nu}} \\
=t^{n(\la)-m\abs{\mu}}
u_0^{(n)}(P_{\mu})\,
\frac{(qt^{m-1};q,t)_{\la}}{c'_{\la}(q,t)}
\prod_{i=1}^m \prod_{j=1}^n 
\frac{(q^{1-N}t^{j-i+m-1}/z;q)_{\la_i-\mu_j}}{(q^{1-N}t^{j-i+m}/z;q)_{\la_i-\mu_j}},
\end{multline*}
where we have also used 
\begin{equation*}
\frac{(a;q)_{N-k}}{(b;q)_{N-k}}=
\frac{(a;q)_N}{(b;q)_N}\,
\frac{(q^{1-N}/b;q)_k}{(q^{1-N}/a;q)_k}\Bigl(\frac{b}{a}\Bigr)^k
\end{equation*}
to rewrite the double product on the right.

Specialising $z\to q^{-N}t^n$ eliminates all reference to $N$ and
completes the proof.
\end{proof}

\section{$\A{n}$ basic hypergeometric series}\label{sec3}

In this section we will be working with $n$ different sets of
variables $\xa{1},\dots,\xa{n}$ where $$\xs=(\xs_1,\dots,\xs_{k_s})$$
such that $k_1\leq k_2\leq \dots \leq k_n$.

Our main object of interest is the following generalisation
of the Kaneko--Mac\-do\-nald basic hypergeometric series \eqref{Phirs}.
\begin{definition}[$\A{n}$ basic hypergeometric series]
\begin{align}\label{seriesdef}
&{_{r+1}\Phi_r}\biggl[\genfrac{}{}{0pt}{}{a_1,\dots,a_{r+1}}
{b_1,\dots,b_r};q,t;\xa{1},\dots,\xa{n}\biggr] \\
&\qquad=\sum_{\laa{1},\dots,\laa{n}} 
\frac{(a_1,\dots,a_{r+1};q,t)_{\laa{n}}}
{(qt^{k_n-1},b_1,\dots, b_r;q,t)_{\laa{n}}} 
\notag \\
&\qquad\qquad\qquad\times
\prod_{s=1}^n 
\biggl[ t^{n(\las)}
\frac{(qt^{k_s-1};q,t)_{\las}}{c'_{\las}(q,t)}\,
P_{\las}(\xs;q,t)\biggr] 
\notag \\
&\qquad\qquad\qquad\times
\prod_{s=1}^{n-1}
\prod_{i=1}^{k_s}\prod_{j=1}^{k_{s+1}}
\frac{(q t^{j-i+k_s-k_{s+1}-1};q)_{\las_i-\laa{s+1}_j}}
{(q t^{j-i+k_s-k_{s+1}};q)_{\las_i-\laa{s+1}_j}}. \notag 
\end{align}
Here the sum is over partitions $\las$ such that
$l(\las)\leq k_s$ for $1\leq s\leq n$ and
\begin{equation}\label{order}
\las_i\geq \laa{s+1}_{i-k_s+k_{s+1}} \quad
\text{for~~$1\leq i\leq k_s$.}
\end{equation}
\end{definition}

\noindent\textbf{Remark.}
The sum over the partitions $\la^{(1)},\dots,\la^{(n)}$ 
subject to \eqref{order} may alternatively be viewed as 
a sum over skew plane partitions of shape $\eta-\nu$ with
$\eta=(k_n^n)$ a partition of rectangular shape and
$\nu=(k_n-k_1,k_n-k_2,\dots,k_n-k_{n-1})$.

\medskip

The above definition simplifies to \eqref{Phirs} when $n=1$, and to
\begin{multline*}
{_{r+1}\Phi_r}\biggl[\genfrac{}{}{0pt}{}{a_1,\dots,a_{r+1}}
{b_1,\dots,b_r};q,t;(z_1),(z_2),\dots,(z_n)\biggr] \\
=\sum_{0\leq j_n\leq \dots\leq j_1}
\frac{(a_1,\dots,a_{r+1};q)_{j_n}}
{(q,b_1,\dots, b_r;q)_{j_n}} \,
z_1^{j_1}\cdots z_n^{j_n}
\prod_{s=1}^{n-1}
\frac{(q/t;q)_{j_s-j_{s+1}}} {(q;q)_{j_s-j_{s+1}}}
\end{multline*}
when $k_1=k_2=\dots=k_n=1$.
Introducing new summation indices by $m_s=j_s-j_{s+1}$ for $1\leq s\leq n$
(where $j_{n+1}:=0$) this gives
\begin{multline*}
{_{r+1}\Phi_r}\biggl[\genfrac{}{}{0pt}{}{a_1,\dots,a_{r+1}}
{b_1,\dots,b_r};q,t;(z_1),(z_2),\dots,(z_n)\biggr] \\
={_{r+1}\phi_r}\biggl[\genfrac{}{}{0pt}{}{a_1,\dots,a_{r+1}}
{b_1,\dots,b_r};q,z_1\cdots z_n\biggr]
\prod_{s=1}^{n-1}
{_1\phi_0}\biggl[\genfrac{}{}{0pt}{}{q/t}{\text{--}}\,;q,z_1\cdots z_s\biggr].
\end{multline*}
Summing the $_1\phi_0$ series by the
$q$-binomial theorem \cite[Equation (II.3)]{GR04}
\begin{equation*}
{_1\phi_0}\biggl[\genfrac{}{}{0pt}{}{a}{\text{--}}\,;q,z\biggr]
=\frac{(az;q)_{\infty}}{(z;q)_{\infty}}
\end{equation*}
we get
\begin{multline}\label{een}
{_{r+1}\Phi_r}\biggl[\genfrac{}{}{0pt}{}{a_1,\dots,a_{r+1}}
{b_1,\dots,b_r};q,t;(z_1),(z_2),\dots,(z_n)\biggr] \\
={_{r+1}\phi_r}\biggl[\genfrac{}{}{0pt}{}{a_1,\dots,a_{r+1}}
{b_1,\dots,b_r};q,z_1\cdots z_n\biggr]
\prod_{s=1}^{n-1}\frac{(qz_1\cdots z_s/t;q)_{\infty}}
{(z_1\cdots z_s;q)_{\infty}}.
\end{multline}
Using Theorem~\ref{thm1} this may be generalised as follows.
\begin{theorem}\label{thm2}
Let $\xa{1}=(\xa{1}_1,\dots,\xa{1}_{k_1})$
and 
\begin{equation*}
\xs=z_s(1,t,\dots,t^{k_s-1})\quad\text{~~ for $2\leq s\leq n$.}
\end{equation*}
Then
\begin{multline*}
{_{r+1}\Phi_r}\biggl[\genfrac{}{}{0pt}{}{a_1,\dots,a_{r+1}}
{b_1,\dots,b_r};q,t;\xa{1},\dots,\xa{n}\biggr] \\
={_{r+1}\Phi_r}\biggl[\genfrac{}{}{0pt}{}{a_1,\dots,a_{r+1}}
{b_1,\dots,b_r};q,t;\hat{x}^{(n)}\biggr]
\prod_{s=1}^{n-1}
\prod_{i=1}^{k_s}\frac{(q\hat{x}^{(s)}_it^{k_s-k_{s+1}-1};q)_{\infty}}
{(\hat{x}^{(s)}_i;q)_{\infty}}.
\end{multline*}
Here the $\hat{x}^{(s)}$ are recursively defined as
$\hat{x}^{(1)}=\xa{1}$ and
\begin{equation*}
\hat{x}^{(s)}=z_s (t^{k_{s-1}-1} \hat{x}^{(s-1)},t^{k_{s-1}},t^{k_{s-1}+1},
\dots,t^{k_s-1})
\quad\text{for~~$2\leq s\leq n$}. 
\end{equation*}
\end{theorem}
Taking $k_1=k_2=\cdots=k_n=1$ and $\xa{1}=(z_1)$ (so that 
$\hat{x}^{(s)}=(z_1\cdots z_s)$) this reduces to \eqref{een}.

Before presenting a proof we will give several important consequences of 
Theorem~\ref{thm2}.
\begin{theorem}[$\A{n}$ $q$-binomial theorem]\label{thm3}
With the same notation as in Theorem~\ref{thm2}
\begin{equation*}
{_1\Phi_0}\biggl[\genfrac{}{}{0pt}{}{a}{\text{--}}\,;q,t;\xa{1},\dots,\xa{n}\biggr]
=\prod_{i=1}^{k_n}\frac{(a\hat{x}_i^{(n)};q)_{\infty}}
{(\hat{x}_i^{(n)};q)_{\infty}}
\prod_{s=1}^{n-1}
\prod_{i=1}^{k_s}\frac{(q\hat{x}^{(s)}_it^{k_s-k_{s+1}-1};q)_{\infty}}
{(\hat{x}^{(s)}_i;q)_{\infty}}.
\end{equation*}
\end{theorem}
Eliminating the hats from the double product on the right yields 
\begin{align*}
{_1\Phi_0}&\biggl[\genfrac{}{}{0pt}{}{a}{\text{--}}\,;q,t;\xa{1},\dots,\xa{n}\biggr] \\
&=\prod_{i=1}^{k_1}\biggl[
\frac{(az_2\cdots z_n\xa{1}_it^{k_1+\cdots+k_{n-1}-n+1};q)_{\infty}}
{(z_2\cdots z_n\xa{1}_it^{k_1+\cdots+k_{n-1}-n+1};q)_{\infty}} \\
&\qquad\qquad\times\prod_{r=1}^{n-1}
\frac{(qz_2\cdots z_r\xa{1}_it^{k_1+\cdots+k_r-k_{r+1}-r};q)_{\infty}}
{(z_2\cdots z_r\xa{1}_it^{k_1+\cdots+k_{r-1}-r+1};q)_{\infty}}
\biggr] \\
&\quad\times
\prod_{s=2}^n \prod_{i=1}^{k_{s}-k_{s-1}}
\frac{(az_s\cdots z_nt^{i+s+k_{s-1}+\cdots+k_{n-1}-n-1};q)_{\infty}}
{(z_s\cdots z_nt^{i+s+k_{s-1}+\cdots+k_{n-1}-n-1};q)_{\infty}} \\
&\quad\times
\prod_{2\leq s\leq r\leq n-1} \prod_{i=1}^{k_s-k_{s-1}}
\frac{(qz_s\cdots z_rt^{i+s-r+k_{s-1}+\cdots+k_r-k_{r+1}-2};q)_{\infty}}
{(z_s\cdots z_rt^{i+s-r+k_{s-1}+\cdots+k_{r-1}-1};q)_{\infty}}.
\end{align*}

\begin{proof}[Proof of Theorem~\ref{thm3}]
If we take Theorem~\ref{thm2} with $r=1$ the ${_1\Phi_0}$
on the right may be summed by \eqref{Phi10sum}, leading to
the desired result.
\end{proof}

If we further specialise $\xa{1}=z_1(1,t,\dots,t^{k_1})$ in Theorem~\ref{thm3}
we obtain a more symmetric $q$-binomial theorem.
\begin{corollary}\label{cor1}
Let $\xs=z_s(1,t,\dots,t^{k_s-1})$ for $1\leq s\leq n$ and set $k_0=0$.
Then
\begin{multline}\label{cor1eq}
{_1\Phi_0}\biggl[\genfrac{}{}{0pt}{}{a}{\text{--}}\,;q,t;\xa{1},\dots,\xa{n}\biggr] \\
=\prod_{s=1}^n \prod_{i=1}^{k_s-k_{s-1}}
\frac{(az_s\cdots z_nt^{i+s+k_{s-1}+\cdots+k_{n-1}-n-1};q)_{\infty}}
{(z_s\cdots z_nt^{i+s+k_{s-1}+\cdots+k_{n-1}-n-1};q)_{\infty}} \\
\times
\prod_{1\leq s\leq r\leq n-1}\prod_{i=1}^{k_s-k_{s-1}}
\frac{(qz_s\cdots z_rt^{i+s-r+k_{s-1}+\cdots+k_r-k_{r+1}-2};q)_{\infty}}
{(z_s\cdots z_rt^{i+s-r+k_{s-1}+\cdots+k_{r-1}-1};q)_{\infty}}.
\end{multline}
\end{corollary}
When $k_1=\dots=k_n=k$ the above significantly simplifies to
\begin{multline*}
{_1\Phi_0}\biggl[\genfrac{}{}{0pt}{}{a}{\text{--}\,};q,t;\xa{1},\dots,\xa{n}\biggr] \\
=\prod_{i=1}^k\biggl[
\frac{(az_1\cdots z_nt^{i-1+(n-1)(k-1)};q)_{\infty}}
{(z_1\cdots z_nt^{i-1+(n-1)(k-1)};q)_{\infty}}
\prod_{s=1}^{n-1} \frac{(qz_1\cdots z_st^{i-2+(s-1)(k-1)};q)_{\infty}}
{(z_1\cdots z_st^{i-1+(s-1)(k-1)};q)_{\infty}}\biggr],
\end{multline*}
where $\xs=z_s(1,t,\dots,t^{k-1})$ for $1\leq s\leq n$.

\noindent\textbf{Remark.}
It is again easily seen that Theorem~\ref{thm3} and Corollary~\ref{cor1}
are true as functions of $\xa{1},z_2,\dots,z_n$ or
$z_1,z_2,\dots,z_n$ when $\abs{q}<1$ and
$$\max\{\abs{\xa{1}_1},\dots,\abs{\xa{1}_{k_1}},
\abs{z_2},\dots,\abs{z_n}\}<1$$
or
$$\max\{\abs{z_1},\dots,\abs{z_n}\}<1.$$

\begin{proof}[Proof of Theorem~\ref{thm2}]
We abbreviate the sequences
$a_1,\dots,a_{r+1}$ and $b_1,\dots,b_r$
by $A$ and $B$, and assume that $n>1$.

If we apply identity \eqref{mn} to eliminate the double product
$\prod_{i=1}^{k_s}\prod_{j=1}^{k_{s+1}}$
in the definition \eqref{seriesdef} of the $_{r+1}\Phi_r$ series
we obtain
\begin{align}\label{nsum}
&{_{r+1}\Phi_r}\biggl[\genfrac{}{}{0pt}{}{A}{B};q,t;\xa{1},\dots,\xa{n}\biggr]\\
&\qquad=\sum_{\substack{\laa{1},\dots,\laa{n}\\
\nu^{(1)},\dots,\nu^{(n-1)} \\
\omega^{(1)},\dots,\omega^{(n-1)}}}
t^{n(\laa{n})} 
\frac{P_{\laa{n}}(\xa{n};q,t)}{c'_{\laa{n}}(q,t)} \,
\frac{(A;q,t)_{\laa{n}}}{(B;q,t)_{\laa{n}}} \notag \\
&\qquad\qquad\qquad\times
\prod_{s=1}^{n-1} 
\biggl[ t^{n(\nu^{(s)})+k_s\abs{\laa{s+1}}-\abs{\omega^{(s)}}}
f_{\omega^{(s)}\nu^{(s)}}^{\las}(q,t)
P_{\las}(\xs;q,t) \notag  \\
&\qquad\qquad\qquad\qquad\quad\times
\frac{u_0^{(k_{s+1}-k_s)}(P_{\laa{s+1}/\omega^{(s)}})}
{u_0^{(k_{s+1})}(P_{\laa{s+1}})}\,
\frac{(qt^{k_s-k_{s+1}-1};q,t)_{\nu^{(s)}}}{c'_{\nu^{(s)}}(q,t)}
\biggr]. \notag 
\end{align}
The $\laa{1}$-dependent part of the summand is given by
\begin{equation*}
f_{\omega^{(1)}\nu^{(1)}}^{\laa{1}}(q,t) P_{\laa{1}}(\xa{1};q,t).
\end{equation*}
Hence the sum over $\laa{1}$ may be performed by \eqref{qtLR} 
to yield
\begin{align*}
&{_{r+1}\Phi_r}\biggl[\genfrac{}{}{0pt}{}{A}{B};q,t;\xa{1},\dots,\xa{n}\biggr]\\
&\qquad=\sum_{\substack{\laa{2},\dots,\laa{n}\\
\nu^{(2)},\dots,\nu^{(n-1)} \\
\omega^{(2)},\dots,\omega^{(n-1)}}}
t^{n(\laa{n})} \frac{P_{\laa{n}}(\xa{n};q,t)}{c'_{\laa{n}}(q,t)}\,
\frac{(A;q,t)_{\laa{n}}}{(B;q,t)_{\laa{n}}}  \,
\frac{t^{k_1\abs{\laa{2}}}}{u_0^{(k_2)}(P_{\laa{2}})} \\ 
&\qquad\qquad\qquad\times
\prod_{s=2}^{n-1}
\biggl[ t^{n(\nu^{(s)})+{k_s}\abs{\laa{s+1}}-\abs{\omega^{(s)}}}
f_{\omega^{(s)}\nu^{(s)}}^{\las}(q,t)
P_{\las}(\xs;q,t) \\
&\qquad\qquad\qquad\qquad\qquad\times
\frac{u_0^{(k_{s+1}-k_s)}(P_{\laa{s+1}/\omega^{(s)}})}
{u_0^{(k_{s+1})}(P_{\la{s+1}})}\,
\frac{(qt^{k_s-k_{s+1}-1};q,t)_{\nu^{(s)}}}{c'_{\nu^{(s)}}(q,t)}
\biggr] \\
&\qquad\qquad\qquad\times
\sum_{\omega^{(1)}}
t^{-\abs{\omega^{(1)}}}
u_0^{(k_2-k_1)}(P_{\laa{2}/\omega^{(1)}})
P_{\omega^{(1)}}(\xa{1};q,t) \\
&\qquad\qquad\qquad\times
\sum_{\nu^{(1)}}
t^{n(\nu^{(1)})}
\frac{(qt^{k_1-k_2-1};q,t)_{\nu^{(1)}}}{c'_{\nu^{(1)}}(q,t)}\,
P_{\nu^{(1)}}(\xa{1};q,t).
\end{align*}
By \eqref{hom} and \eqref{skewdef2} the sum over $\omega^{(1)}$ gives
\begin{equation*}
P_{\laa{2}}(t^{-1}\xa{1},1,\dots,t^{k_2-k_1-1};q,t)
\end{equation*}
and by \eqref{Phirs} and \eqref{Phi10sum}
the sum over $\nu^{(1)}$ gives
\begin{equation*}
{_1\Phi_0}\biggl[\genfrac{}{}{0pt}{}
{qt^{k_1-k_2-1}}{\text{--}}\,;q,t;\xa{1}\biggr]
=\prod_{i=1}^{k_1}\frac{(q\xa{1}_it^{k_1-k_2-1};q)_{\infty}}
{(\xa{1}_i;q)_{\infty}}.
\end{equation*}
Substituting these two results and once again using
\eqref{hom} to absorb the factor $t^{k_1\abs{\laa{2}}}$, we find
\begin{align*}
&{_{r+1}\Phi_r}\biggl[\genfrac{}{}{0pt}{}{A}{B};q,t;\xa{1},\dots,\xa{n}\biggr]
=\prod_{i=1}^{k_1}\frac{(q\xa{1}_it^{k_1-k_2-1};q)_{\infty}}{(\xa{1}_i;q)_{\infty}} \\
&\qquad\times
\sum_{\substack{\laa{2},\dots,\laa{n}\\
\nu^{(2)},\dots,\nu^{(n-1)} \\
\omega^{(2)},\dots,\omega^{(n-1)}}}
t^{n(\laa{n})}
\frac{P_{\laa{n}}(\xa{n};q,t)}{c'_{\laa{n}}(q,t)}\,
\frac{(A;q,t)_{\laa{n}}}{(B;q,t)_{\laa{n}}} \\
&\qquad\qquad\qquad\times
\frac{P_{\laa{2}}(t^{k_1-1}\xa{1},t^{k_1},\dots,t^{k_2-1};q,t)}
{u_0^{(k_2)}(P_{\laa{2}})} \\
&\qquad\qquad\qquad\times
\prod_{s=2}^{n-1}
\biggl[ t^{n(\nu^{(s)})+k_s\abs{\laa{s+1}}-\abs{\omega^{(s)}}}
f_{\omega^{(s)}\nu^{(s)}}^{\las}(q,t) 
P_{\las}(\xs;q,t) \\
&\qquad\qquad\qquad\qquad\qquad\times
\frac{u_0^{(k_{s+1}-k_s)}(P_{\laa{s+1}/\omega^{(s)}})}
{u_0^{(k_{s+1})}(P_{\laa{s+1}})}\,
\frac{(qt^{k_s-k_{s+1}-1};q,t)_{\nu^{(s)}}}{c'_{\nu^{(s)}}(q,t)}
\biggr].
\end{align*}
Comparing this with \eqref{nsum} we see that up to the term
\begin{equation}\label{term}
\frac{P_{\laa{2}}(t^{k_1-1}\xa{1},t^{k_1},\dots,t^{k_2-1};q,t)}
{u_0^{(k_2)}(P_{\laa{2}})}
\end{equation}
we have effectively reduced $n$ to $n-1$.
The naive approach would be to apply $u_{0;t^{1-k_1}}^{(k_1)}$ acting on $\xa{1}$.
Then \eqref{term} collapses to $1$ as desired, but 
\begin{equation*}
u_{0;t^{1-k_1}}^{(k_1)}\biggl(\,
\prod_{i=1}^{k_1}\frac{(\xa{1}_it^{k_1-k_2-1};q)_{\infty}}
{(\xa{1}_i;q)_{\infty}}\biggr)=
\prod_{i=1}^{k_1}\frac{(t^{k_1-k_2-i};q)_{\infty}}
{(t^{1-i};q)_{\infty}}
\end{equation*}
is not well-defined.
The correct approach is to apply $u_{0;z_2}^{(k_2)}$ acting on $\xa{2}$,
resulting in
\begin{align*}
&u_{0;z_2}^{(k_2)}\biggl({_{r+1}\Phi_r}\biggl[\genfrac{}{}{0pt}{}{A}{B};q,t;
\xa{1},\xa{2},\xa{3},\dots,\xa{n}\biggr] \biggr)\\
&\quad=
\prod_{i=1}^{k_1}\frac{(q\xa{1}_it^{k_1-k_2-1};q)_{\infty}}{(\xa{1}_i;q)_{\infty}} \\
&\qquad\times
\sum_{\substack{\laa{2},\dots,\laa{n}\\
\nu^{(2)},\dots,\nu^{(n-1)} \\
\omega^{(2)},\dots,\omega^{(n-1)}}}
t^{n(\laa{n})}
\frac{P_{\laa{n}}(\xa{n};q,t)}{c'_{\laa{n}}(q,t)}\,
\frac{(A;q,t)_{\laa{n}}}{(B;q,t)_{\laa{n}}} \\
&\qquad\qquad\qquad\times
\prod_{s=2}^{n-1}
\biggl[ t^{n(\nu^{(s)})+k_s\abs{\laa{s+1}}-\abs{\omega^{(s)}}}
f_{\omega^{(s)}\nu^{(s)}}^{\las}(q,t)
P_{\las}(\hat{x}^{(s)};q,t) \\
&\qquad\qquad\qquad\qquad\qquad\times
\frac{u_0^{(k_{s+1}-k_s)}(P_{\laa{s+1}/\omega^{(s)}})}
{u_0^{(k_{s+1})}(P_{\laa{s+1}})}\,
\frac{(qt^{k_s-k_{s+1}-1};q,t)_{\nu^{(s)}}}{c'_{\nu^{(s)}}(q,t)}
\biggr],
\end{align*}
where $\hat{x}^{(2)}=z_2(t^{k_1-1}\xa{1},t^{k_1},\dots,t^{k_2-1})$
and $\hat{x}^{(s)}=\xs$ for $3\leq s\leq n-1$.
Again comparing this with \eqref{nsum} we have thus proved the following
intermediate result.
\begin{lemma}
Assume that $n>1$ and let $\xa{2}=z_2(1,t,\dots,t^{k_2-1})$
and $\hat{x}^{(2)}=z_2(t^{k_1-1}\xa{1},t^{k_1},\dots,t^{k_2-1})$.
Then
\begin{multline*}
{_{r+1}\Phi_r}\biggl[\genfrac{}{}{0pt}{}{A}{B};q,t;\xa{1},\dots,\xa{n}\biggr] \\
={_{r+1}\Phi_r}\biggl[\genfrac{}{}{0pt}{}{A}{B};q,t;\hat{x}^{(2)},\xa{3},\dots,\xa{n}\biggr]
\prod_{i=1}^{k_1}\frac{(q\xa{1}_it^{k_1-k_2-1};q)_{\infty}}{(\xa{1}_i;q)_{\infty}}.
\end{multline*}
\end{lemma}
This is readily iterated, resulting in Theorem~\ref{thm2}.
\end{proof}

\section{Proof of the $\A{n}$ exponential Selberg integral}\label{sec4}

Although Corollary~\ref{thmESelbergAn} follows as a straightforward limit
of the $\A{n}$ Selberg integral, it is proved here
directly from Corollary~\ref{cor1}.
The advantage of first dealing with the exponential integral
instead of the full $\A{n}$ Selberg integral is that
it makes the introduction of 
$C^{k_1,\dots,k_n}_{\gamma}[a,b]$ slightly simpler.
Throughout the proof we assume that $\gamma\neq 0$.

Applying $u_{0;z_1}^{(k_1)}\cdots u_{0;z_n}^{(k_n)}$ (with 
$u_{0;z_s}^{(k_s)}$ acting on $\xs$) to \eqref{seriesdef},
and using \eqref{ccp} and \eqref{PS}, we get
\begin{align}\label{som}
&{_{r+1}\Phi_r}\biggl[\genfrac{}{}{0pt}{}{A}{B};q,t;
\xa{1},\dots,\xa{n}\biggr]
=\sum_{\laa{1},\dots,\laa{n}} 
\frac{(A;q,t)_{\laa{n}}}
{(qt^{k_n-1},B;q,t)_{\laa{n}}} \\
&\quad\times
\prod_{s=1}^n 
\biggl[ t^{2n(\las)}z_s^{\abs{\las}}
\prod_{1\leq i<j\leq k_s}
\frac{1-t^{j-i}q^{\las_i-\las_j}}{1-t^{j-i}}
\frac{(t^{j-i+1};q)_{\las_i-\las_j}}{(qt^{j-i-1};q)_{\las_i-\las_j}} \biggr] 
\notag \\
&\quad\times
\prod_{s=1}^{n-1}
\prod_{i=1}^{k_s}\prod_{j=1}^{k_{s+1}}
\frac{(q t^{j-i+k_s-k_{s+1}-1};q)_{\las_i-\laa{s+1}_j}}
{(q t^{j-i+k_s-k_{s+1}};q)_{\las_i-\laa{s+1}_j}},
\notag 
\end{align}
where $\xs=z_s(1,t,\dots,t^{k_s-1})$ for $1\leq s\leq n$.
Taking $r=0$ and $A=a_1=a$ this may be equated with the right-hand side of
\eqref{cor1eq}. Then replacing $t\to q^{\gamma}$ and
$a\to q^{\alpha+(k_n-1)\gamma}$, and letting $q\to 1^{-}$ using
\begin{equation*}
(q^x;q)_k \to \frac{\Gamma(x+k)}{\Gamma(x)}
\end{equation*}
and
\begin{equation}\label{zyx}
\frac{(q^xz;q)_{\infty}}{(q^yz;q)_{\infty}} \to 
(1-z)^{y-x}  \qquad \abs{z}<1,
\end{equation}
yields
\begin{align}\label{q1sum}
&\sum_{\laa{1},\dots,\laa{n}} 
\prod_{i=1}^{k_n}
\frac{\Gamma(\alpha+\muu{n}_i)}{\Gamma(1+\muu{n}_i)} 
\prod_{s=1}^{n-1}\prod_{i=1}^{k_s}\prod_{j=1}^{k_{s+1}}
\frac{\Gamma(1-\gamma+\mus_i-\muu{s+1}_j)}
{\Gamma(1+\mus_i-\muu{s+1}_j)}  \\
&\qquad\qquad\times
\prod_{s=1}^n 
\biggl[z_s^{\abs{\las}}
\prod_{1\leq i<j\leq k_s}
\frac{(\mus_i-\mus_j)\,\Gamma(\gamma+\mus_i-\mus_j)}
{\Gamma(1-\gamma+\mus_i-\mus_j)} \biggr]\notag \\
&\quad =
\prod_{1\leq s\leq r\leq n}
(1-z_s\cdots z_r)^{-(\alpha_r+(k_r-k_{r+1}-1)\gamma)(k_s-k_{s-1})} 
\notag \\
&\quad\quad\times
\prod_{s=1}^n \prod_{i=1}^{k_s}
\frac{\Gamma(\alpha_s+(i-k_{s+1}-1)\gamma)\Gamma(i\gamma)}{\Gamma(\gamma)}.
\notag 
\end{align}
Here $\mus_i:=\las_i+(k_s-i)\gamma$, $k_0=k_{n+1}=0$ and
$\alpha_s$ is defined as in Theorem~\ref{thmSelbergAn}.

Next we would like to replace
$z_s\to\exp(-\epsilon\beta_s)$ and $\la_i^{(s)}\to t_i^{(s)}/\epsilon$
(with $\epsilon>0$),
and take the $\epsilon\to 0$ limit using Stirling's formula
\begin{equation*}
\lim_{x\to\infty} x^{b-a} \frac{\Gamma(x+a)}{\Gamma(x+b)}=1
\end{equation*}
to transform the above sum into an integral.
There is however the complication that
the difference $\las_i-\laa{s+1}_j$
is not necessarily nonnegative. Let 
\begin{equation}\label{udef}
u_{ij}^{(s)}=i-j-k_s+k_{s+1}.
\end{equation}
Then for $\las_i-\laa{s+1}_j>0$ 
\begin{multline*}
\frac{\Gamma(1-\gamma+\mus_i-\muu{s+1}_j)}
{\Gamma(1+\mus_i-\muu{s+1}_j)} \\
\to
\frac{\Gamma(1-(u_{ij}^{(s)}+1)\gamma+(t_i^{(s)}-t_j^{(s+1)})/\epsilon)}
{\Gamma(1-u_{ij}^{(s)}\gamma+(t_i^{(s)}-t_j^{(s+1)})/\epsilon)}
\sim \epsilon^{\gamma}\bigl(t_i^{(s)}-t_j^{(s+1)}\bigr)^{-\gamma},
\end{multline*}
but for $\las_i-\laa{s+1}_j<0$, by the Euler reflection formula,
\begin{align*}
\frac{\Gamma(1-\gamma+\mus_i-\muu{s+1}_j)}
{\Gamma(1+\mus_i-\muu{s+1}_j)}
&=\frac{\Gamma(-\mus_i+\muu{s+1}_j)}
{\Gamma(\gamma-\mus_i+\muu{s+1}_j)} \,
R_{ij}^{(s)}(\gamma) \\
&\to
\frac{\Gamma(u_{ij}^{(s)}\gamma-(t_i^{(s)}-t_j^{(s+1)})/\epsilon)}
{\Gamma((u_{ij}^{(s)}+1)\gamma-(t_i^{(s)}-t_j^{(s+1)})/\epsilon)} \,
R_{ij}^{(s)}(\gamma) \\
& \sim \epsilon^{\gamma}
\bigl(t_j^{(s+1)}-t_i^{(s)}\bigr)^{-\gamma} R_{ij}^{(s)}(\gamma) \\
& =\epsilon^{\gamma}
\abs{t_i^{(s)}-t_j^{(s+1)}}^{-\gamma} R_{ij}^{(s)}(\gamma),
\end{align*}
where 
\begin{equation}\label{R}
R_{ij}^{(s)}(\gamma)=
\frac{\sin\bigl(\pi u_{ij}^{(s)}\gamma\bigr)}
{\sin\bigl(\pi(u_{ij}^{(s)}+1)\gamma\bigr)}.
\end{equation}
Hence, 
depending on the relative order of $\las_i$ and $\laa{s+1}_j$ we 
may or may not pick up a ratio of sine-functions, and
for small $\epsilon$ the summand of \eqref{q1sum} takes the form
\begin{multline*}
\epsilon^{k_1+\cdots+k_n-\sum_{1\leq s\leq r \leq n}
(\alpha_r+(k_r-k_{r+1}-1)\gamma)(k_s-k_{s-1})}  
\prod R_{ij}^{(s)}(\gamma)
\\ \times
\prod_{i=1}^{k_n}
(t_i^{(n)})^{\alpha-1}
\prod_{s=1}^{n-1}
\abs{\Delta(t^{(s)},t^{(s+1)})}^{-\gamma} 
\prod_{s=1}^n\eup^{-\beta_s\sum_{i=1}^{k_s}t_i^{(s)}} 
\bigl(\Delta(t^{(s)})\bigr)^{2\gamma},
\end{multline*}
where the product is over all $i,j$ and $s$ such that
$t_i^{(s)}<t_j^{(s+1)}$.
{}From this it follows that we must first fix a complete ordering
between the parts of $\las$ and $\laa{s+1}$.
Any such ordering compatible with \eqref{order} may be described 
by a map \cite{TV03}
\begin{subequations}\label{Map}
\begin{equation}
M_s:\{1,\dots,k_s\}\to \{1,\dots,k_{s+1}\}
\end{equation}
such that 
\begin{equation}\label{C1}
M_s(i)\leq M_s(i+1)
\end{equation}
and
\begin{equation}\label{C2}
1\leq M_s(i)\leq k_{s+1}-k_s+i,
\end{equation}
\end{subequations}
as follows:
\begin{equation}\label{ineq}
\laa{s+1}_{M_s(i)}\leq \las_i\leq \laa{s+1}_{M_s(i)-1}
\quad \text{for~~$1\leq i\leq k_s$},
\end{equation}
where $\laa{s+1}_0=\infty$.
Note in particular that \eqref{C1} must hold so that
\eqref{ineq} is compatible with the ordering among the
parts of $\las$. Similarly \eqref{C2} must hold so that
\eqref{ineq} is compatible with \eqref{order}.
A straightforward counting argument shows that there are
exactly
\begin{equation*}
\frac{k_{s+1}-k_s+1}{k_{s+1}+1}\binom{k_s+k_{s+1}}{k_s}
\end{equation*}
different maps.

Now define $D^{k_1,\dots,k_n}[a,b]\subseteq [a,b]^{k_1+\cdots+k_n}$ 
as the set of points
\begin{equation*}
(\laa{1},\dots,\laa{n})=
(\laa{1}_1,\dots,\laa{1}_{k_1},\dots,\laa{n}_1,\dots,\laa{n}_{k_n})
\end{equation*}
satisfying
\begin{equation*}
a\leq \las_{k_s}\leq \dots\leq\las_1\leq b
\quad \text{for~~$1\leq s\leq n$}
\end{equation*}
and \eqref{order}.
Given admissible maps $M_1,\dots,M_{n-1}$ we define
$D_{M_1,\dots,M_{n-1}}^{k_1,\dots,k_n}[a,b]\subseteq
D^{k_1,\dots,k_n}[a,b]$ by requiring that
\eqref{ineq} holds for $1\leq s\leq n-1$.
As chains we have
\begin{equation}\label{domain}
D^{k_1,\dots,k_n}[a,b]=\sum_{M_1,\dots,M_{n-1}}
D_{M_1,\dots,M_{n-1}}^{k_1,\dots,k_n}[a,b],
\end{equation}
and summing over $\la^{(1)},\dots,\la^{(n)}$ amounts to
summing over the 
lattice points in
$D^{k_1,\dots,k_n}[0,\infty]$.
Thanks to the decomposition \eqref{domain} we know exactly
which factors of the form $R_{ij}^{(s)}(\gamma)$ are 
picked up when we go from sum to integral.
Indeed, from \eqref{ineq} we have that
$\las_i\leq \laa{s+1}_j$ for $1\leq j\leq M_s(i)-1$.
This gives rise to the factor
\begin{equation*}
\prod_{j=1}^{M_s(i)-1}
R_{ij}^{(s)}(\gamma)
=\frac{\sin\bigl(\pi(i+k_{s+1}-k_s-M_s(i)+1)\gamma\bigr)}
{\sin\bigl(\pi(i+k_{s+1}-k_s)\gamma\bigr)}.
\end{equation*}
Taking the product over $i$ and $s$ this yields
the total factor
\begin{equation*}
F_{M_1,\dots,M_{n-1}}^{k_1,\dots,k_n}(\gamma)=
\prod_{s=1}^{n-1}\prod_{i=1}^{k_s}
\frac{\sin\bigl(\pi(i+k_{s+1}-k_s-M_s(i)+1)\gamma\bigr)}
{\sin\bigl(\pi(i+k_{s+1}-k_s)\gamma\bigr)}.
\end{equation*}
Hence making the variable changes
$z_s\to \exp(-\epsilon\beta_s)$ and $\la_i^{(s)}\to t_i^{(s)}/\epsilon$
and letting $\epsilon$ tend to zero we obtain Corollary~\ref{thmESelbergAn},
where the integration is over the chain
\begin{equation}\label{chain}
C^{k_1,\dots,k_n}_{\gamma}[a,b]=
\sum_{M_1,\dots,M_{n-1}}
F_{M_1,\dots,M_{n-1}}^{k_1,\dots,k_n}(\gamma)
D_{M_1,\dots,M_{n-1}}^{k_1,\dots,k_n}[a,b].
\end{equation}
For $n=2$ this corresponds to the chain introduced in
by Tarasov and Varchenko in \cite{TV03}
(up to some trivial notational changes).
Of course, to correctly interpret \eqref{chain} in the context of 
Corollary~\ref{thmESelbergAn} (and also Theorem~\ref{thmSelbergAn})
we have to replace $\las_i$ by $t_i^{(s)}$ in all of the above.
In particular \eqref{ineq} becomes
\begin{equation*}
t^{(s+1)}_{M_s(i)}\leq t^{(s)}_i\leq t^{(s+1)}_{M_s(i)-1}
\quad \text{for~~$1\leq i\leq k_s$},
\end{equation*}
and we integrate the $t_i^{(s)}$ such that
\begin{equation*}
(t_1^{(1)},\dots,t_{k_1}^{(1)},\dots,t_1^{(n)},\dots,t_{k_n}^{(n)})
\in D^{k_1,\dots,k_n}[a,b].
\end{equation*}

\section{Proof of the $\A{n}$ Selberg integral}\label{sec5}
Throughout this section we assume $\gamma\neq 0$ and $0<q<1$, and use
\[
(a;q)_z=\frac{(a;q)_{\infty}}{(aq^z;q)_{\infty}}
\quad \text{for $z\in\Complex$}.
\]

Let $\Gamma_q$ be the $q$-gamma function \cite[Equation (I.35)]{GR04}
\begin{equation*}
\Gamma_q(x)=\frac{(q;q)_{x-1}}{(1-q)^{x-1}},
\end{equation*}
and define
\begin{equation*}
\Delta_{\gamma}\bigl(\xs;q\bigr)=\prod_{1\leq i<j\leq k_s}
(\xs_j)^{2\gamma}\bigl(1-q^{(j-i)\gamma}\xs_i/\xs_j\bigr)
\bigl(q^{1+(j-i-1)\gamma}\xs_i/\xs_j;q\bigr)_{2\gamma-1}
\end{equation*}
and
\begin{equation*}
\Delta_{\gamma}\bigl(\xs,\xa{s+1};q\bigr)=
\prod_{i=1}^{k_s}\prod_{j=1}^{k_{s+1}}
\bigl(\xa{s+1}_j\bigr)^{-\gamma}
\bigl(q^{1-u_{ij}^{(s)}\gamma}\xs_i/\xa{s+1}_j;q\bigr)_{-\gamma},
\end{equation*}
with $u_{ij}^{(s)}$ as in \eqref{udef}.
\begin{comment}
Then 
\begin{align*}
\lim_{q\to 1^{-}}\Gamma_q(x)&=\Gamma(x), \\
\lim_{q\to 1^{-}} \Delta_{\gamma}(\xs;q)&=
\prod_{1\leq i<j\leq k_s}(\xs_j-\xs_i)^{2\gamma}
\intertext{and}
\lim_{q\to 1^{-}}\Delta_{\gamma}(\xs,\xa{s+1};q)&=
\prod_{i=1}^{k_s}\prod_{j=1}^{k_{s+1}}(x_j^{(s+1)}-x_i^{(s)})^{-\gamma}.
\end{align*}
\end{comment}

Using the above definitions as well as \eqref{som},
Corollary~\ref{cor1} can be written in the form
\begin{align}\label{qi}
&(1-q)^{k_1+\cdots+k_n}\sum_{\laa{1},\dots,\laa{n}}
\prod_{i=1}^{k_n} \bigl(q^{1+(k_n-i)\gamma}\xa{n}_i;q\bigr)_{\alpha-1}
\prod_{s=1}^n\prod_{i=1}^{k_s}\bigl(\xs_i\bigr)^{\beta_s} \\
&\qquad\qquad\qquad\qquad\qquad
\times \prod_{s=1}^{n-1} \Delta_{\gamma}\bigl(\xs,\xa{s+1};q\bigr)
\prod_{s=1}^n \Delta_{\gamma}\bigl(\xs;q\bigr)\notag \\
&=\prod_{1\leq s\leq r\leq n} \prod_{i=1}^{k_s-k_{s-1}}
\frac{\Gamma_q(\beta_s+\cdots+\beta_r+(i+s-r-1)\gamma)}
{\Gamma_q(\alpha_r+\beta_s+\cdots+\beta_r+(i+s-r+k_r-k_{r+1}-2)\gamma)}
\notag \\
&\quad \times
\prod_{s=1}^n \prod_{i=1}^{k_s}
\frac{\Gamma_q(\alpha_s+(i-k_{s+1}-1)\gamma)\Gamma_q(i\gamma)}
{\Gamma_q(\gamma)}. \notag
\end{align}
Here $\alpha_s$ is as defined in Theorem~\ref{thmSelbergAn}, $\xs_i:=q^{\las_i}$
and
$(t,z_s,a)$ has been replaced by $(q^{\gamma},q^{\beta_s-k_{s-1}\gamma},
q^{(k_n-1)\gamma+\alpha})$.

The reader familiar with Jackson or $q$-integrals will recognise \eqref{qi} 
as a $(k_1+\cdots+k_n)$-dimensional such integral.
The standard $1$-dimensional Jackson integral is given by
\begin{equation*}
\int_0^1 f(x)\dup_q x=(1-q) \sum_{i=0}^{\infty} f(q^i) q^i
\end{equation*}
which simplifies (at least formally) to the Riemann integral
\begin{equation*}
\int_0^1 f(x)\dup x
\end{equation*}
in the $q\to 1^{-}$ limit.
Generalising this to arbitrary dimensions as
\begin{equation*}
\Int_{[0,1]^n} f(x)\dup_q x=(1-q)^n \sum_{i_1,\dots,i_n=0}^{\infty} 
f(q^{i_1},\dots,q^{i_n}) q^{i_1+\cdots+i_n},
\end{equation*}
where $x=(x_1,\dots,x_n)$ and $\dup_q x=\dup_q x_1\cdots \dup_q x_n$,
\eqref{qi} corresponds to the restricted $q$-integral
\begin{align}\label{Dint}
&\Int_{\bar{D}^{k_1,\dots,k_n}[0,1]} 
\prod_{i=1}^{k_n} \bigl(q^{1+(k_n-i)\gamma}\xa{n}_i;q\bigr)_{\alpha-1}
\prod_{s=1}^n\prod_{i=1}^{k_s}\bigl(\xs_i\bigr)^{\beta_s-1} \\
&\qquad\qquad\qquad
\times \prod_{s=1}^{n-1} \Delta_{\gamma}\bigl(\xs,\xa{s+1};q\bigr)
\prod_{s=1}^n \Delta_{\gamma}\bigl(\xs;q\bigr) \:
\dup_q \xa{1} \cdots \dup_q \xa{n} \notag \\
&=\prod_{1\leq s\leq r\leq n} \prod_{i=1}^{k_s-k_{s-1}}
\frac{\Gamma_q(\beta_s+\cdots+\beta_r+(i+s-r-1)\gamma)}
{\Gamma_q(\alpha_r+\beta_s+\cdots+\beta_r+(i+s-r+k_r-k_{r+1}-2)\gamma)}
\notag  \\
&\quad \times
\prod_{s=1}^n \prod_{i=1}^{k_s}
\frac{\Gamma_q(\alpha_s+(i-k_{s+1}-1)\gamma)\Gamma_q(i\gamma)}
{\Gamma_q(\gamma)}, \notag
\end{align}
where 
$\bar{D}^{k_1,\dots,k_n}[a,b]\subseteq [a,b]^{k_1+\cdots+k_n}$ is the set of all points
\begin{equation*}
\bigl(\xa{1},\dots,\xa{n}\bigr)=
\bigl(\xa{1}_1,\dots,\xa{1}_{k_1},\dots,\xa{n}_1,\dots,\xa{n}_{k_n}\bigr)
\end{equation*}
such that
\begin{equation*}
a\leq \xs_1\leq \cdots\leq \xs_{k_s}\leq b \quad 
\text{for~~$1\leq s\leq n$}
\end{equation*}
and
\begin{equation*}
\xs_i\leq \xa{s+1}_{i-k_s+k_{s+1}} \quad
\text{for~~$1\leq i\leq k_s$,~~$1\leq s\leq n-1$.}
\end{equation*}

Assuming that
\begin{subequations}\label{con}
\begin{equation}
\xs_i<\xs_j \quad \text{for~~$1\leq i<j\leq k_s$}
\end{equation}
and
\begin{equation}\label{conb}
\xs_i<\xa{s+1}_j \quad
\text{for~~$1\leq i\leq k_s$,~~$1\leq j\leq k_{s+1}$}
\end{equation}
\end{subequations}
it follows from \eqref{zyx} that
\begin{subequations}
\begin{align}
\lim_{q\to 1^{-}} \Delta_{\gamma}\bigl(\xs;q\bigr)&=
\bigl(\Delta\bigl(-\xs\bigr)\bigr)^{2\gamma} , \\
\lim_{q\to 1^{-}}\Delta_{\gamma}\bigl(\xs,\xa{s+1};q\bigr)&=
\bigl(\Delta\bigl(-\xs,-\xa{s+1}\bigr)\bigr)^{-\gamma}.
\label{limb}
\end{align}
\end{subequations}
Hence, for \eqref{con} (and $\xa{n}_i<1$ for $1\leq i\leq k_n$)
the $q\to 1^{-}$ limit of the integrand of \eqref{Dint} is
\begin{equation}\label{integrand}
\prod_{i=1}^{k_n} \bigl(1-\xa{n}_i\bigr)^{\alpha-1}
\prod_{s=1}^n\biggl[\,
\bigl(\Delta\bigl(-\xs\bigr)\bigr)^{2\gamma}
\prod_{i=1}^{k_s}\bigl(\xs_i\bigr)^{\beta_s-1}\, \biggr]
\prod_{s=1}^{n-1} \bigl(\Delta\bigl(-\xs,-\xa{s+1}\bigr)\bigr)^{-\gamma}.
\end{equation}

The problem we are now facing is exactly the same as that of
Section~\ref{sec4}: the $\xs_i$ are not necessarily 
compatible with \eqref{conb}.
This forces us to also consider the limit of the factors making up 
$\Delta_{\gamma}(\xs,\xa{s+1};q)$ when $\xs_i>\xa{s+1}_j$.
In computing the limit \eqref{limb} we used that
\begin{equation*}
\lim_{q\to 1^{-}}
\bigl(q^{1-u_{ij}^{(s)}\gamma}\xs_i/\xa{s+1}_j;q\bigr)_{-\gamma}
=\bigl(1-\xs_i/\xa{s+1}_j\bigr)^{-\gamma},
\end{equation*}
but this is only correct for $\xs_i<\xa{s+1}_j$.
When $\xs_i>\xa{s+1}_j$ we may use
the $q$-reflection formula \cite[Equation (168a)]{Thomae}
\[
\Gamma_q(z)\Gamma_q(1-z)=\frac{2\iup q^{z/2} \theta_1(i\log q^{z/2};q^{1/2})}
{(1-q) \theta_1'(0;q^{1/2})}
\]
(upon recalling that $\xa{s}_i:=q^{\las_i}$) to write
\[
\bigl(q^{1-u_{ij}^{(s)}\gamma}\xs_i/\xa{s+1}_j;q\bigr)_{-\gamma}=
\frac{\bigl(q^{-1/2}\xa{s+1}_j/\xa{s}_i\bigr)^{\gamma}}
{\bigl(q^{u_{ij}^{(s)}\gamma}\xa{s+1}_j/\xs_i;q\bigr)_{\gamma}} 
\;\frac{\theta_1(\iup \log q^{u_{ij}^{(s)}\gamma/2};q^{1/2})}
{\theta_1(\iup \log q^{(u_{ij}^{(s)}+1)\gamma/2};q^{1/2})}.
\]
Since
\begin{multline*}
\frac{q^{u/2}\,\theta_1(\iup \log q^{u/2};q^{1/2})}
{q^{v/2}\,\theta_1(\iup \log q^{v/2};q^{1/2})}
=\frac{(1-q^u)(q^{1+u};q)_{\infty}(q^{1-u};q)_{\infty}}
{(1-q^v)(q^{1+v};q)_{\infty}(q^{1-v};q)_{\infty}} \\
\stackrel{q\to 1^{-}}{\longrightarrow} \:
\frac{u}{v}\prod_{n=1}^{\infty}\frac{1-u^2/n^2}{1-v^2/n^2}
=\frac{\sin\pi u}{\sin\pi v}
\end{multline*}
we find
\begin{align*}
\lim_{q\to 1^{-}}
\bigl(q^{1-u_{ij}^{(s)}\gamma}\xs_i/\xa{s+1}_j;q\bigr)_{-\gamma}
&=\bigl(\xs_i/\xa{s+1}_j-1\bigr)^{-\gamma} R_{ij}^{(s)}(\gamma) \\
&=\Abs{1-\xs_i/\xa{s+1}_j}^{-\gamma} R_{ij}^{(s)}(\gamma),
\end{align*}
where $R_{ij}^{(s)}$ is given by \eqref{R}.
We therefore conclude that the $q\to 1^{-}$ limit of the integrand of
\eqref{Dint} is given by \eqref{integrand} with 
\begin{equation*}
\bigl(\Delta\bigl(-\xs,-\xa{s+1}\bigr)\bigr)^{-\gamma}
\to
\Abs{\Delta\bigl(-\xs,-\xa{s+1}\bigr)}^{-\gamma}
=\Abs{\Delta\bigl(\xs,\xa{s+1}\bigr)}^{-\gamma}.
\end{equation*}
multiplied by a factor $R_{ij}^{(s)}(\gamma)$
for each $\xs_i>\xa{s+1}_j$.

The rest of the proof simply follows Section~\ref{sec4}.
We introduce the maps $M_s$ as in \eqref{Map},
with $M_s(i)$ such that (compare with \eqref{ineq})
\begin{equation}\label{ineq2}
\xa{s+1}_{M_s(i)-1}\leq \xs_i\leq \xa{s+1}_{M_s(i)}
\quad \text{for~~$1\leq i\leq k_s$},
\end{equation}
where $\xa{s+1}_0=0$.
Furthermore we define $\bar{D}^{k_1,\dots,k_n}_{M_1,\dots,M_n}[a,b]\subseteq
\bar{D}^{k_1+\cdots+k_n}[a,b]$ by requiring that \eqref{ineq2} holds. Therefore, 
\begin{equation}\label{domain2}
\bar{D}^{k_1,\dots,k_n}[a,b]=\sum_{M_1,\dots,M_{n-1}}
\bar{D}_{M_1,\dots,M_{n-1}}^{k_1,\dots,k_n}[a,b]
\end{equation}
as chains.
Also defining 
\begin{equation}\label{chain2}
\bar{C}^{k_1,\dots,k_n}_{\gamma}[a,b]=
\sum_{M_1,\dots,M_{n-1}}
F_{M_1,\dots,M_{n-1}}^{k_1,\dots,k_n}(\gamma)
\bar{D}_{M_1,\dots,M_{n-1}}^{k_1,\dots,k_n}[a,b]
\end{equation}
we get
\begin{align*}
&\Int_{\bar{C}^{k_1,\dots,k_n}_{\gamma}[0,1]}
\prod_{i=1}^{k_n}
\bigl(1-x_i^{(n)}\bigr)^{\alpha-1}
\prod_{s=1}^n\biggl[\,
\Abs{\Delta\bigl(x^{(s)}\bigr)}^{2\gamma}
\prod_{i=1}^{k_s}\bigl(x_i^{(s)}\bigr)^{\beta_s-1}\,\biggr] \\
&\qquad\qquad\quad\times \prod_{s=1}^{n-1}
\Abs{\Delta\bigl(x^{(s)},x^{(s+1)}\bigr)}^{-\gamma} 
\; \dup \xa{1}\cdots\dup \xa{n} \\
&=\prod_{1\leq s\leq r\leq n} \prod_{i=1}^{k_s-k_{s-1}}
\frac{\Gamma(\beta_s+\cdots+\beta_r+(i+s-r-1)\gamma)}
{\Gamma(\alpha_r
+\beta_s+\cdots+\beta_r+(i+s-r+k_r-k_{r+1}-2)\gamma)} \\
&\quad\times
\prod_{s=1}^n \prod_{i=1}^{k_s}
\frac{\Gamma(\alpha_s+(i-k_{s+1}-1)\gamma)
\Gamma(i\gamma)}{\Gamma(\gamma)}
\end{align*}
(where $\abs{\Delta(x^{(s)})}^{2\gamma}$ may be replaced by
$(\Delta(-x^{(s)}))^{2\gamma}$\,).
Finally making the variable changes $\xs_i=1-t^{(s)}_i$ for all $1\leq i\leq k_s$ 
$1\leq s\leq n$, so that
$\bar{C}^{k_1,\dots,k_n}_{\gamma}[0,1]$ is replaced by
$C^{k_1,\dots,k_n}_{\gamma}[0,1]$, completes the proof.

\section{Further $\A{n}$ integrals}\label{sec6}

Let $P_{\la}^{(\alpha)}$ be the Jack polynomial, obtained from
the Macdonald polynomial $P_{\la}$ as
\begin{equation*}
P_{\la}^{(\alpha)}(x)=\lim_{t\to 1} P_{\la}(x;t^{\alpha},t),
\end{equation*}
and let $(a)_N$ 
be the Pochhammer symbol
\[
(a)_N=a(a+1)\cdots(a+N-1).
\]
Then it is an easy matter to generalise the previous derivation to yield an
$\A{n}$ integral involving the Jack polynomial.
\begin{theorem}\label{thmJack}
Let $\mu$ be a partition of at most $k_1$ parts.
With the same conditions as in Theorem~\ref{thmSelbergAn} we have
\begin{align*}
&\Int_{\bar{C}^{k_1,\dots,k_n}_{\gamma}[0,1]}
P_{\mu}^{(1/\gamma)}\bigl(\xa{1}\bigr)
\prod_{s=1}^n  \biggl[\,
\Abs{\Delta\bigl(x^{(s)}\bigr)}^{2\gamma}
\prod_{i=1}^{k_s}\bigl(1-x_i^{(s)}\bigr)^{\alpha_s-1}
\bigl(x_i^{(s)}\bigr)^{\beta_s-1}\,\biggr] \\
&\qquad\qquad\quad\times
\prod_{s=1}^{n-1}
\Abs{\Delta\bigl(x^{(s)},x^{(s+1)}\bigr)}^{-\gamma}
\; \dup \xa{1}\cdots\dup \xa{n} \\
&=\prod_{1\leq i<j\leq k_1}\frac{((j-i+1)\gamma)_{\mu_i-\mu_j}}
{((j-i)\gamma)_{\mu_i-\mu_j}} \\
&\quad\times
\prod_{s=1}^n
\prod_{i=1}^{k_1} 
\frac{(\beta_1+\cdots+\beta_s+(k_1-s-i+1)\gamma)_{\mu_i}}
{(\alpha_s+\beta_1+\cdots+\beta_s+(k_1+k_s-k_{s+1}-s-i)\gamma)_{\mu_i}} \\
&\quad \times
\prod_{1\leq s\leq r\leq n} \prod_{i=1}^{k_s-k_{s-1}}
\frac{\Gamma(\beta_s+\cdots+\beta_r+(i+s-r-1)\gamma)}
{\Gamma(\alpha_r
+\beta_s+\cdots+\beta_r+(i+s-r+k_r-k_{r+1}-2)\gamma)} \\
&\quad\times
\prod_{s=1}^n \prod_{i=1}^{k_s}
\frac{\Gamma(\alpha_s+(i-k_{s+1}-1)\gamma)
\Gamma(i\gamma)}{\Gamma(\gamma)}.
\end{align*}
\end{theorem}
When $n=1$ this simplifies to
\begin{align}\label{kadell}
&\frac{1}{k!}\Int_{[0,1]^k}P_{\mu}^{(1/\gamma)}(x)
\abs{\Delta(x)}^{2\gamma}
\prod_{i=1}^k x_i^{\alpha-1}(1-x_i)^{\beta-1} \; \dup x \\
&\qquad=\prod_{1\leq i<j\leq k}\frac{\Gamma((j-i+1)\gamma+\mu_i-\mu_j)}
{\Gamma((j-i)\gamma+\mu_i-\mu_j)} \notag \\
&\quad\qquad \times
\prod_{i=1}^k
\frac{\Gamma(\alpha+(k-i)\gamma+\mu_i)\Gamma(\beta+(i-1)\gamma)}
{\Gamma(\alpha+\beta+(2k-i-1)\gamma+\mu_i)}, \notag
\end{align}
where we have made the substitutions
$(k_1,\alpha,\beta_1)\to (k,\beta,\alpha)$ and have used the
symmetry of the integrand to replace
\begin{equation*}
\Int_{\bar{C}^k_{\gamma}[0,1]}=\Int_{0\leq x_1\leq x_2\dots\leq x_k\leq 1}
\quad\text{by}\quad\quad \frac{1}{k!}\Int_{[0,1]^k}.
\end{equation*}
The integral \eqref{kadell} is due to Kadell \cite[Theorem I]{Kadell97}
(see also \cite[pp. 385--386]{Macdonald95}).
The special case $\mu=(1^r)$ of \eqref{kadell} corresponds to Aomoto's
integral \cite{Aomoto87}, usually stated as
\begin{multline*}
\Int_{[0,1]^k}
\abs{\Delta(x)}^{2\gamma}
\prod_{i=1}^r x_i 
\prod_{i=1}^k x_i^{\alpha-1} (1-x_i)^{\beta-1}\, \dup x  \\
=\prod_{i=1}^r \frac{(\alpha+(k-i)\gamma)}{(\alpha+\beta+(2k-i-1)\gamma)}
\prod_{i=1}^k
\frac{\Gamma(\alpha+(i-1)\gamma) \Gamma(\beta+(i-1)\gamma)\Gamma(i\gamma+1)}
{\Gamma(\alpha+\beta+(i+k-2)\gamma)\Gamma(\gamma+1)},
\end{multline*}
for $0\leq r\leq k$.
The equivalence of Aomoto's integral and the $\mu=(1^r)$ case of
\eqref{kadell} follows by symmetrising the integrand of the former, 
using that
\begin{equation}\label{elem}
\sum_{i_1<i_2<\dots<i_r} x_{i_1}x_{i_2}\cdots x_{i_r}=
e_r(x)=P^{(1/\gamma)}_{(1^r)}(x),
\end{equation}
with $e_r$ the $r$th elementary symmetric function.

By taking $\mu=(1^r)$ in Theorem~\ref{thmJack} we obtain the 
following $\A{n}$ analogue of Aomoto's integral:
\begin{align*}
&\Int_{\bar{C}^{k_1,\dots,k_n}_{\gamma}[0,1]}e_r(\xa{1})
\prod_{s=1}^n
\biggl[\,
\Abs{\Delta\bigl(x^{(s)}\bigr)}^{2\gamma}
\prod_{i=1}^{k_s} \bigl(1-x_i^{(s)}\bigr)^{\alpha_s-1}
\bigl(x_i^{(s)}\bigr)^{\beta_s-1}\,\biggr] \\
&\qquad\qquad\quad\times
\prod_{s=1}^{n-1}
\Abs{\Delta\bigl(x^{(s)},x^{(s+1)}\bigr)}^{-\gamma}
\; \dup \xa{1}\cdots\dup \xa{n} \\
&=\binom{k_1}{r}
\prod_{s=1}^n\prod_{i=1}^r\frac{(\beta_1+\cdots+\beta_s+(k_1-i-s+1)\gamma)}
{(\alpha_s+\beta_1+\cdots+\beta_s+(k_1+k_s-k_{s+1}-i-s)\gamma)} \\
&\quad \times
\prod_{1\leq s\leq r\leq n} \prod_{i=1}^{k_s-k_{s-1}}
\frac{\Gamma(\beta_s+\cdots+\beta_r+(i+s-r-1)\gamma)}
{\Gamma(\alpha_r
+\beta_s+\cdots+\beta_r+(i+s-r+k_r-k_{r+1}-2)\gamma)} \\
&\quad\times
\prod_{s=1}^n \prod_{i=1}^{k_s}
\frac{\Gamma(\alpha_s+(i-k_{s+1}-1)\gamma)
\Gamma(i\gamma)}{\Gamma(\gamma)}.
\end{align*}

\begin{proof}[Proof of Theorem~\ref{thmJack}]
In proving Theorem~\ref{thmSelbergAn} we have not taken advantage of the
full $\A{n}$ $q$-binomial theorem as stated in Theorem~\ref{thm3},
relying on the less general Corollary~\ref{cor1} instead.
In going from the former to the latter we have specialised 
$\xa{1}$ to $z_1(1,t,\dots,t^{k_1-1})$, or, equivalently,
applied to $u_{0;z_1}^{(k_1)}$ acting on $\xa{1}$.

If, more generally, we apply $u_{\mu;z_1}^{(k_1)}$
instead of $u_{0;z_1}^{(k_1)}$, the
factor $u_0^{(k_1)}(P_{\laa{1}})$ in the summand is replaced
by $u_{\mu}^{(k_1)}(P_{\laa{1}})$.
Then invoking \eqref{symm} this leads to the term
\begin{equation*}
u_0^{(k_1)}(P_{\laa{1}})\: 
\frac{u_{\laa{1}}^{(k_1)}(P_{\mu})}{u_0^{(k_1)}(P_{\mu})}
\end{equation*}
instead of just $u_0^{(k_1)}(P_{\laa{1}})$.

Of course not just the summand of \eqref{cor1eq} will change by the
above, and by applying $u_{\mu;z_1}^{(k_1)}$ instead of $u_{0;z_1}^{(k_1)}$,
the right-hand side of \eqref{cor1eq} picks up the additional factor
\begin{equation*}
\prod_{s=1}^n
\frac{(z_1\cdots z_st^{k_1+\cdots+k_{s-1}+k_1-s};q,t)_{\mu}}
{(a_sz_1\cdots z_st^{k_1+\cdots+k_s-k_{s+1}+k_1-s-1};q,t)_{\mu}},
\end{equation*}
where $a_1=\cdots=a_{n-1}=q$ and $a_n=at^{1-k_n}$.
Accordingly, the identity \eqref{qi} generalises to
\begin{align*}
&(1-q)^{k_1+\cdots+k_n}\sum_{\laa{1},\dots,\laa{n}}
P_{\mu}(y;q,q^{\gamma})
\prod_{i=1}^{k_n} \bigl(q^{1+(k_n-i)\gamma}\xa{n}_i;q\bigr)_{\alpha-1}
\prod_{s=1}^n\prod_{i=1}^{k_s}\bigl(\xs_i\bigr)^{\beta_s} \\
&\qquad\qquad\qquad\qquad\qquad
\times \prod_{s=1}^{n-1} \Delta_{\gamma}\bigl(\xs,\xa{s+1};q\bigr)
\prod_{s=1}^n \Delta_{\gamma}\bigl(\xs;q\bigr)\notag \\
&=u_0^{(k_1)}(P_{\mu}(q,q^{\gamma}))
\prod_{s=1}^n\frac{(q^{\beta_1+\cdots+\beta_s+(k_1-s)\gamma};
q,q^{\gamma})_{\mu}}
{(q^{\alpha_s+\beta_1+\cdots+\beta_s+(k_1+k_s-k_{s+1}-s-1)\gamma};
q,q^{\gamma})_{\mu}} \\
&\quad \times
\prod_{1\leq s\leq r\leq n} \prod_{i=1}^{k_s-k_{s-1}}
\frac{\Gamma_q(\beta_s+\cdots+\beta_r+(i+s-r-1)\gamma)}
{\Gamma_q(\alpha_r+\beta_s+\cdots+\beta_r+(i+s-r+k_r-k_{r+1}-2)\gamma)}
\\
&\quad \times
\prod_{s=1}^n \prod_{i=1}^{k_s}
\frac{\Gamma_q(\alpha_s+(i-k_{s+1}-1)\gamma)\Gamma_q(i\gamma)}
{\Gamma_q(\gamma)},
\end{align*}
where $y=(y_1,\dots,y_{k_1})$ and $y_i=\xa{1}_i q^{(k_i-i)\gamma}$.
The rest of the proof proceeds exactly as before.
\end{proof}

\section{Two simple examples}\label{sec7}
To end this paper we present the fully worked-out examples
of the $\A{n}$ Selberg integral for 
\[
(k_1,\dots,k_{n-1},k_n)=(1,\dots,1,k)
\]
and for $\gamma=0$.
\subsection{The case $(k_1,\dots,k_{n-1},k_n)=(1,\dots,1,k)$}

In this case there is only one map $M_s$ for $1\leq s\leq n-2$, 
corresponding to the identity map $M_s(1)=1$.
For $s=n-1$, however, there are $k$ different maps, 
given by $M_{n-1}(1)=a$ for $1\leq a\leq k$.

If we relabel the integration variables $t_1^{(s)}\to u_s$ 
for $1\leq s\leq n-1$ and $t_i^{(n)}\to t_i$ for $1\leq i\leq k$, 
then the above implies the inequalities
\begin{equation*}
O_a:
\begin{cases}
0\leq t_k\leq\dots\leq t_1\leq 1,\\
0\leq u_{n-1}\leq \dots\leq u_1\leq 1,\\
t_a\leq u_{n-1}\leq t_{a-1}
\end{cases}
\end{equation*}
with $1\leq a\leq k$ and $t_0=1$.
As a result we obtain the following $(k+n-1)$-dimensional integral:
\begin{align}\label{example}
&\sum_{a=1}^k 
\frac{\sin(\pi(k-a+1)\gamma)}{\sin(\pi k\gamma)} \\
&\quad\times
\int_{O_a}\:
\prod_{i=1}^{n-1} (1-u_i)^{\beta_i-1} 
\prod_{i=1}^k
t_i^{\alpha-1} (1-t_i)^{\beta_n-1} 
\prod_{1\leq i<j\leq k}(t_i-t_j)^{2\gamma} \notag \\
&\qquad\quad\times
\prod_{i=1}^{n-2}(u_i-u_{i+1})^{-\gamma}
\prod_{i=1}^{a-1}(t_i-u_{n-1})^{-\gamma}
\prod_{i=a}^k(u_{n-1}-t_i)^{-\gamma} 
\: \dup u \, \dup t  \notag \\
&=\Gamma(1-k\gamma)\Gamma^{n-2}(1-\gamma)
\prod_{i=1}^k \frac{\Gamma(\alpha+(i-1)\gamma)\Gamma(i\gamma)}
{\Gamma(\gamma)} \notag \\
&\quad\times
\prod_{i=1}^{k-1}
\frac{\Gamma(\beta_n+(i-1)\gamma)}
{\Gamma(\alpha+\beta_n+(i+k-2)\gamma)} 
\prod_{i=1}^n\frac{\Gamma(\beta_1+\cdots+\beta_i+(1-i)\gamma)}
{\Gamma(A_i+\beta_1+\cdots+\beta_i-i\gamma)}, \notag 
\end{align}
where $A_1=\cdots=A_{n-2}=1$, $A_{n-1}=1-(k-1)\gamma$,
$A_n=\alpha+k\gamma$, $\dup u=\dup u_1 \cdots \dup u_{n-1}$ and
$\dup t=\dup t_1 \cdots \dup t_k$.

In the notation of the introduction the above integral corresponds to
\begin{equation*}
I_{\underbrace{\scriptstyle 1,\dots,1}_{n-1},k}^{\A{n}}
(\alpha;\beta_1,\dots,\beta_n;\gamma)
\end{equation*}
and, according to the recurrence \eqref{keen},
all but one of the ones may be eliminated.

To see this assume that $n\geq 3$ and replace the integration variable
$u_1$ by $v$ as
\begin{equation*}
v=\frac{u_1-u_2}{1-u_2}.
\end{equation*}
Noting that $1-u_1=(1-v)(1-u_2)$ and $u_1-u_2=v(1-u_2)$ the integral over
$v$ may be identified as Euler's beta integral \eqref{Euler} with 
$\alpha=1-\gamma$.
Therefore
\begin{equation*}
I_{\underbrace{\scriptstyle 1,\dots,1}_{n-1},k}^{\A{n}}
(\alpha;\beta_1,\dots,\beta_n;\gamma)
=I_{\underbrace{\scriptstyle 1,\dots,1}_{n-2},k}^{\A{n-1}}
(\alpha;\beta_1+\beta_2-\gamma,\beta_3,\dots,\beta_n;\gamma)\,
\frac{\Gamma(1-\gamma)\Gamma(\beta_1)}{\Gamma(\beta_1-\gamma+1)}.
\end{equation*}
in accordance with \eqref{keen}.
Iterating the recursion it follows that
\begin{multline*}
I_{\underbrace{\scriptstyle 1,\dots,1}_{n-1},k}^{\A{n}}
(\alpha;\beta_1,\dots,\beta_n;\gamma)
=I_{1,k}^{\A{2}}(\alpha;\beta_1+\cdots+\beta_n-(n-1)\gamma,\beta_n;\gamma)  \\
\times \Gamma^{n-2}(1-\gamma)\prod_{i=1}^{n-2}
\frac{\Gamma(\beta_1+\cdots+\beta_i+(1-i)\gamma)}
{\Gamma(A_i+\beta_1+\cdots+\beta_i-i\gamma)}
\end{multline*}
and \eqref{example} boils down to its $\A{2}$ or $n=2$ case
\begin{align*}
&\sum_{a=1}^k 
\frac{\sin(\pi(k-a+1)\gamma)}{\sin(\pi k\gamma)} \\
&\quad\times
\int_{O_a'}\:
(1-u)^{\beta_1-1} 
\prod_{i=1}^k
t_i^{\alpha-1} (1-t_i)^{\beta_2-1} 
\prod_{1\leq i<j\leq k}(t_i-t_j)^{2\gamma} \\
&\qquad\quad\times
\prod_{i=1}^{a-1}(t_i-u)^{-\gamma} \prod_{i=a}^k(u-t_i)^{-\gamma} 
\: \dup u \, \dup t  \\
&=\frac{\Gamma(\beta_1)\Gamma(1-k\gamma)}{\Gamma(1+\beta_1-k\gamma)}\,
\frac{\Gamma(\alpha+\beta_2+(2k-2)\gamma)}{\Gamma(\alpha+\beta_1+\beta_2+(k-2)\gamma)}\,
\frac{\Gamma(\beta_1+\beta_2-\gamma)}{\Gamma(\beta_2+(k-1)\gamma)} \\
&\quad\times
\prod_{i=1}^k \frac{\Gamma(\alpha+(i-1)\gamma)\Gamma(\beta_2+(i-1)\gamma)
\Gamma(i\gamma)}
{\Gamma(\alpha+\beta_2+(i+k-2)\gamma)\Gamma(\gamma)} .
\end{align*}
where 
\begin{equation*}
O_a':\quad 0\leq t_k\leq\dots\leq t_a\leq u\leq t_{a-1}\leq\dots\leq t_1\leq 1.
\end{equation*}
\begin{comment}
($k=1$ is trivial) except for the degenerate case $\gamma=0$
corresponding to
\begin{multline*}
\sum_{a=1}^k \frac{k-a+1}{k} 
\int_{O_a'}\: (1-u)^{\beta_1-1} \prod_{i=1}^k
t_i^{\alpha-1} (1-t_i)^{\beta_2-1} \: \dup u \, \dup t  \\
=\frac{1}{\beta_1 k!}\,
\frac{\Gamma(\alpha)\Gamma(\beta_1+\beta_2)}{\Gamma(\alpha+\beta_1+\beta_2)}\,
\biggl(\frac{\Gamma(\alpha)\Gamma(\beta_2)}{\Gamma(\alpha+\beta_2)}\biggr)^{k-1} .
\end{multline*}
After integrating over $u$ this gives
\begin{multline*}
\sum_{a=1}^k \frac{(k-a+1)}{k}
\Int_{0\leq t_k\leq\dots\leq t_1\leq 1}
\Bigl[(1-t_a)^{\beta_1}-
(1-t_{a-1})^{\beta_1}\Bigr] 
\prod_{i=1}^k
t_i^{\alpha-1} (1-t_i)^{\beta_2-1} \: \dup t \\
=\frac{1}{k!}
\frac{\Gamma(\alpha)\Gamma(\beta_1+\beta_2)}{\Gamma(\alpha+\beta_1+\beta_2)}\,
\biggl(\frac{\Gamma(\alpha)\Gamma(\beta_2)}{\Gamma(\alpha+\beta_2)}\biggr)^{k-1} .
\end{multline*}
where $t_0=1$. Denoting the left-hand side by $\text{LHS}$ 
the derivation is completed by noting that
\begin{align*}
\text{LHS}&=\frac{1}{k}
\Int_{0\leq t_k\leq\dots\leq t_1\leq 1}
\biggl[\:\sum_{a=1}^k (1-t_a)^{\beta_1}\biggr] 
\prod_{i=1}^k
t_i^{\alpha-1} (1-t_i)^{\beta_2-1} \: \dup t \\
&=\frac{1}{k\, k!}
\Int_{[0,1]^k}
\biggl[\:\sum_{a=1}^k (1-t_a)^{\beta_1}\biggr] 
\prod_{i=1}^k
t_i^{\alpha-1} (1-t_i)^{\beta_2-1} \: \dup t \\
&=\frac{1}{k!}
\Int_{[0,1]^k} 
t_1^{\alpha-1} (1-t_1)^{\beta_1+\beta_2-1}
\prod_{i=2}^k
t_i^{\alpha-1} (1-t_i)^{\beta_2-1} \: \dup t
\end{align*}
so that we obtain a product of $k$ Euler beta integrals \eqref{Euler}.
\end{comment}

\subsection{The case $\gamma=0$}
When $\gamma=0$ Theorem~\ref{thm1} collapses to the integral
\begin{multline*}
\Int_{C^{k_1,\dots,k_n}_0[0,1]}
\prod_{s=1}^n \prod_{i=1}^{k_s} 
\bigl(t_i^{(s)}\bigr)^{\alpha_s-1}
\bigl(1-t_i^{(s)}\bigr)^{\beta_s-1} 
 \; \dup t^{(1)}\cdots\dup t^{(n)} \\
=\prod_{s=1}^n \frac{1}{(k_s)!}
\biggl(\frac{\Gamma(\alpha)\Gamma(\beta_s+\cdots+\beta_n)}
{\Gamma(\alpha+\beta_s+\cdots+\beta_n)}\biggr)^{k_s-k_{s-1}}
\end{multline*}
with $\alpha_1=\dots=\alpha_{n-1}=1$ and $\alpha_n=\alpha$.
Because the $t_i^{(s)}$ in the integrand are completely decoupled
the problem of evaluating this integral is purely combinatorial.
Introducing the partitions $\la^{(s)}$ for
$1\leq s\leq n-1$ as $\la^{(s)}=(M_s(k_s),\dots,M_s(1))$
so that $\la^{(s)}$ has exactly $k_s$ parts and $\la_i^{(s)}\leq k_{s+1}-i+1$
the $\gamma=0$ integral may also be stated more explicitly as
\begin{multline*}
\sum_{\substack{\la^{(1)},\dots,\la^{(n-1)} \\
1\leq\la_i^{(s)}\leq k_{s+1}-i+1}}
\prod_{s=1}^{n-1}\prod_{i=1}^{k_s}
\frac{k_{s+1}-i-\la_i+2}{k_{s+1}-i+1} \\ \times
\Int \prod_{s=1}^n \prod_{i=1}^{k_s} 
\bigl(t_i^{(s)}\bigr)^{\alpha_s-1}
\bigl(1-t_i^{(s)}\bigr)^{\beta_s-1} 
\; \dup t^{(1)}\cdots\dup t^{(n)} \\
=\prod_{s=1}^n \frac{1}{(k_s)!}
\biggl(\frac{\Gamma(\alpha)\Gamma(\beta_s+\cdots+\beta_n)}
{\Gamma(\alpha+\beta_s+\cdots+\beta_n)}\biggr)^{k_s-k_{s-1}}
\end{multline*}
where the integration domain is given by
\begin{subequations}\label{IE}
\begin{equation}
\max\Bigl\{t^{(s+1)}_{\la^{(s)}_{k_s-i+1}},t_{i+1}^{(s)}\Bigr\}
\leq t_i^{(s)} \leq 
\min\Bigl\{t_{i-1}^{(s)},t^{(s+1)}_{\la^{(s)}_{k_s-i+1}-1}\Bigr\}
\end{equation}
for $1\leq s\leq n-1$ and $1\leq i\leq k_s$ (with $t_0^{(s)}=1$ and
$t_{k_s+1}^{(s)}=0$), and
\begin{equation}
t_{i+1}^{(n)} \leq t_i^{(n)} \leq t_{i-1}^{(n)}
\end{equation}
\end{subequations}
for $1\leq i\leq k_n$.

Thanks to the factor 
$\prod_{s=1}^{n-1}\prod_{i=1}^{k_s}(k_{s+1}-i-\la_i+2)$ we may relax
the condition $\la_i^{(s)}\leq k_{s+1}-i+1$ to $\la_1^{(s)}\leq k_{s+1}$ 
so that the sum becomes
\[
\sum_{\substack{\la^{(1)},\dots,\la^{(n-1)} \\[1mm]
l(\la^{(s)})=k_s \\[1mm]
\la^{(s)}_1\leq k_{s+1}}}.
\]

Since $\alpha_1=\cdots=\alpha_{n-1}=1$ one may, with a bit of pain,
successively integrate
over the $t_i^{(s)}$ starting with $s=1$. We will not present the
full details of this calculation here, but remark that the key
to unravelling the combinatorics encoded in the inequalities
\eqref{IE} is given by
\begin{equation}\label{er}
\sum_{\substack{\la \\[1mm] l(\la)=r \\[1mm] \la_1\leq n}}
\prod_{j=1}^n \frac{1}{m_j!}
\prod_{i=1}^r (n-\la_i-i+2)(x_{\la_i}-x_{\la_i-1})=
e_r(x),
\end{equation}
with $x_0=0$, $x=(x_1,\dots,x_n)$, 
$m_j$ the multiplicity of the part $j$ in $\lambda$
and $e_r$ the $r$th elementary symmetric function \eqref{elem}.
To establish \eqref{er} we note that when written in terms
of the multiplicities $m_j$ it becomes
\[
\frac{1}{(n-r)!}
\sum_{\substack{m_1,\dots,m_n \geq 0 
\\[1mm] M_1=r}} \:
\prod_{j=1}^n \frac{(x_j-x_{j-1})^{m_j}(n-j+1-M_{j+1})}{m_j!}
=e_r(x),
\]
where $M_j=m_j+\cdots+m_n$. Multiplying this by $t^r$ and summing over 
$r$ using the generating function for the 
$e_r$ \cite[Equation (I.2.2)]{Macdonald95} yields
\[
\sum_{m_1,\dots,m_n \geq 0}\:
\frac{1}{(n-M_1)!}
\prod_{j=1}^n \frac{[t(x_j-x_{j-1})]^{m_j}(n-j+1-M_{j+1})}{m_j!}
=\prod_{i=1}^n(1+tx_i).
\]
This is true for any $x_0$ provided we add the factor
$(1+tx_0)^{n-M_1}$ to the summand:
\begin{equation}\label{x0gen}
\sum_{m_1,\dots,m_n \geq 0}\:
\frac{(1+tx_0)^{n-M_1}}{(n-M_1)!}
\prod_{j=1}^n \frac{[t(x_j-x_{j-1})]^{m_j}(n-j+1-M_{j+1})}{m_j!}
=\prod_{i=1}^n(1+tx_i).
\end{equation}
For $n=0$ this is obviously correct. If we denote the sum on the
left by $L(x_0,\dots,x_n)$ then
\begin{align*}
L(x_0,\dots,x_n)&=
\sum_{m_2,\dots,m_n \geq 0}\:
\frac{(1+tx_0)^{n-M_2}}{(n-1-M_2)!}
\prod_{j=2}^n \frac{[t(x_j-x_{j-1})]^{m_j}(n-j+1-M_{j+1})}{m_j!}\\
&\qquad\qquad\qquad\qquad\times
\sum_{m_1=0}^{n-M_2}\biggl(\frac{tx_1-tx_0}{1+tx_0}\biggr)^{m_1}
\binom{n-M_2}{m_1} \\
&=\sum_{m_2,\dots,m_n \geq 0}\:
\frac{(1+tx_1)^{n-M_2}}{(n-1-M_2)!}
\prod_{j=2}^n \frac{[t(x_j-x_{j-1})]^{m_j}(n-j+1-M_{j+1})}{m_j!}\\
&=(1+tx_1) L(x_1,\dots,x_n).
\end{align*}
By induction \eqref{x0gen} is thus true for all nonnegative integers $n$.

\bibliographystyle{amsplain}

\end{document}